\documentclass[11pt]{article}
\usepackage{latexsym}
\usepackage[all]{xy}
\usepackage{amsmath}
\usepackage{amssymb}
\usepackage{amsfonts}

\newtheorem{thm}{Theorem}[subsection]
\newtheorem{prop}[thm]{Proposition}
\newtheorem{lem}[thm]{Lemma}
\newtheorem{df}[thm]{Definition}
\newtheorem{cor}[thm]{Corollary}
\newtheorem{conj}[thm]{Conjecture}
\newtheorem{hyp}[thm]{Hypothesis}

\begin{document}

\title{\textbf{Algebraic geometry over model categories} \\
A general approach to derived algebraic geometry}
\bigskip
\bigskip

\author{\bigskip\\
Betrand Toen \\
\small{Laboratoire J. A. Dieudonn\'{e}}\\
\small{UMR CNRS 6621} \\
\small{Universit\'{e} de Nice Sophia-Antipolis}\\
\small{France}\\
\bigskip \and
\bigskip \\
Gabriele Vezzosi \\
\small{Dipartimento di Matematica}\\
\small{Universit\`a di Bologna}\\
\small{Italy}\\
\bigskip}

\maketitle

\begin{abstract}
For a (semi-)model category $M$, we define a notion of a
''homotopy'' Grothendieck topology on $M$, as well as its
associated model category of stacks. We use this to define a
notion of geometric stack over a symmetric monoidal base model
category; geometric stacks are the fundamental objects to "do
algebraic geometry over model categories". We give two examples of
applications of this formalism. The first one is the
interpretation of $DG$-schemes as geometric stacks over the model
category of complexes and the second one is a definition of
\'{e}tale $K$-theory of $E_{\infty}$-ring spectra.

This first version is very preliminary and might be considered as
a detailed research announcement. Some proofs, more details and
more examples will be added in a forthcoming version.
\end{abstract}


\tableofcontents

\newpage

\setcounter{section}{0}

\begin{section}{Introduction}

By definition, a scheme is obtained by gluing together affine
schemes for the Zariski topology. Therefore, algebraic geometry is
a theory which is based on the two fundamental notions of
\textit{affine scheme} and \textit{Grothendieck topology}. It was
observed already a long time ago that these two notions still make
sense in more general contexts, and that \textit{schemes} can be
defined in very general settings. This has led to the theory of
\textit{relative algebraic geometry}, which allows one to
\textit{do algebraic geometry} over well behaved symmetric
monoidal base categories (see \cite{del,degal,hak}); usual
algebraic geometry corresponds then to the "absolute" case where
the base category is the category of $\mathbb{Z}$-modules.\\
The goal of the present work is to start a program to develop
\textit{algebraic geometry relatively to symmetric monoidal
$\infty$-categories}. Our motivations for starting such a program
come from several questions in algebraic geometry and algebraic
topology and will be clarified in the two entries \textit{Examples
and applications} and \textit{Relations with other works} of this
introduction.

It is well known that model categories give rise in a natural way
to $\infty$\textit{-categories}. Indeed, B. Dwyer and D. Kan
defined a simplicial localization process, which starting from a
model category $M$, constructs a simplicial category $LM$, the
simplicial localization of $M$ (see \cite{dk1}). As simplicial
categories may be viewed as $\infty$-categories for which
$i$-morphisms are invertible up to $(i+1)$-morphisms for all $i>1$
(see \cite[\S $2$]{sh}), this suggests that model categories
\textit{are} a certain kind of $\infty$-categories. In the same
way, symmetric monoidal model categories (as defined in \cite[\S
$4$]{ho}) \textit{are} a certain kind of symmetric monoidal
$\infty$-categories (e.g. in the sense of \cite{to1}). As a first
step in our program we would like to present in this paper a
setting to \textit{do algebraic geometry relatively to symmetric
monoidal model categories}. For this, we will concentrate on
defining a category of \textit{geometric stacks over a base
symmetric monoidal model category}, whose construction will be the main purpose of this work. \\

\bigskip

\begin{center} \textit{Review of usual algebraic geometry} \end{center}

In order to explain our approach, we shall first present in detail
a construction of the usual category of schemes, or more generally
of algebraic stacks and of $n$-geometric stacks (see \cite{s1}),
emphasizing the categorical ingredients needed in each step, in
such a way that the generalization that will follow should look
fairly natural.

The starting point is the category $Aff$ of \textit{affine
schemes}. By definition, we will take $Aff$ to be the opposite of
the category of commutative and unital rings. We consider its
Yoneda embedding
$$h : Aff \longrightarrow Aff^{\wedge},$$
where $Aff^{\wedge}$ is the category of presheaves of simplicial
sets on $Aff$ (i.e. $Aff^{\wedge}:=\textrm{SPr}(Aff)$) and $h$
maps an affine scheme $X$ to $h_{X}:=Hom(-,X)$ (here a set is
always considered as a constant simplicial set). From a
categorical point of view, the embedding $h : Aff \longrightarrow
Aff^{\wedge}$ is obtained by \textit{formally adding homotopy
colimits to $Aff$} (see \cite{du2} for more details on this point
of view). This process is relevant to our situation, as
\textit{gluing objects in $Aff$} will be done by taking certain
formal homotopy colimits of objects in $Aff$ (i.e. taking homotopy
colimits in $Aff^{\wedge}$ of object in $Aff$). The category
$Aff^{\wedge}$ is in a natural way a model category, where
equivalences are defined objectwise, and the functor $h$ induces a
Yoneda embedding on the level of the homotopy categories
$$h : Aff \longrightarrow Ho(Aff^{\wedge}).$$ Throughout this
introduction, the homotopy category $Ho(\mathcal{C})$ of a
category $\mathcal{C}$ with a distinguished class of morphisms $w$
will denote the category obtained from $\mathcal{C}$ by formally
inverting all morphisms in $w$; when $\mathcal{C}$ is a model
category, we will implicitly assume that $w$ is the set of weak
equivalences and when $\mathcal{C}$ does not come naturally
equipped with a model category structure we will consider it as a
trivial model category where $w$ consists of all the isomorphisms.

The next step is choosing a \textit{Grothendieck topology} on
$Aff$, that will be used to glue affine schemes. For the purpose
of schemes, the Zariski topology is enough, but \'{e}tale or even
faithfully flat and quasi-compact (for short $\textrm{ffqc}$)
topologies also proves very useful in order to define more general
objects as algebraic spaces or algebraic stacks. We will choose
here to work with the $\textrm{ffqc}$ topology though the
construction will be valid for any Grothendieck topology on $Aff$.
The $\textrm{ffqc}$ topology makes $Aff$ into a Grothendieck site
and therefore one can consider its category of
($\infty$-)\textit{stacks}, denoted by
$Ho(Aff^{\sim,\textrm{ffqc}})$. For us, the category
$Ho(Aff^{\sim,\textrm{ffqc}})$ is the full sub-category of
$Ho(Aff^{\wedge})$ consisting of simplicial presheaves satisfying
the descent condition for $\textrm{ffqc}$ hyper-covers
\footnote{For us the word \textit{stack} will always mean a
\textit{stack of $\infty$-groupoids}, adopting the point of view
of \cite{s2} and \cite{to2}, according to which \textit{stacks of
$\infty$-groupoids} are modelled by simplicial presheaves.}. The
category $Ho(Aff^{\sim,\textrm{ffqc}})$ is precisely the homotopy
category of a certain model category structure on $Aff^{\wedge}$,
and the category $Aff^{\wedge}$ together with this model structure
will be called the model category of stacks, and denoted by
$Aff^{\sim,\textrm{ffqc}}$. Finally, it is known that the topology
$\textrm{ffqc}$ is sub-canonical, or in other words that the
Yoneda embedding $h : Aff \longrightarrow Ho(Aff^{\wedge})$
factors through $Ho(Aff^{\sim,\textrm{ffqc}})$. Therefore we have
an induced fully faithful functor
$$h : Aff \longrightarrow Ho(Aff^{\sim,\textrm{ffqc}}).$$
A stack in the essential image of the functor $h$ will be called
by extension an \textit{affine scheme}.

Let us consider now a simplicial object $X_{*} : \Delta^{op}
\longrightarrow Aff^{\sim,\textrm{ffqc}}$ (i.e. $X_{*}$ is a
bi-simplicial presheaf) and suppose that $X_{*}$ satisfies the
following three conditions>

\begin{enumerate}
\item The simplicial object $X_{*}$ is a \textit{Segal groupoid} in $Aff^{\sim,\textrm{\textrm{ffqc}}}$ (see Def.
\ref{d21});

\item The image of each $X_{n}$ in $Ho(Aff^{\sim,ffqc})$ is a disjoint union of affine
schemes;

\item The two morphisms (\textit{source} and \textit{target}) $X_{1} \rightrightarrows X_{0}$ are
faithfully flat and affine morphisms (this makes sense since the
$X_{n}$'s are disjoint union of affine schemes).
\end{enumerate}

To such a groupoid we associate its homotopy colimit $|X_{*}| \in
Ho(Aff^{\sim,\textrm{ffqc}})$, which can be defined to be the
stack associated to the diagonal of the bi-simplicial presheaf
$X_{*}$. It is not difficult to check that the full sub-category
of $Ho(Aff^{\sim,\textrm{ffqc}})$, consisting of objects
isomorphic to some $|X_{*}|$, with $X_{*}$ satisfying conditions
$(1)$, $(2)$ and $(3)$, is equivalent to the homotopy category of
\textit{algebraic stacks} (in the sense of Artin, see \cite{lm})
having an affine diagonal. In particular, it contains the category
of separated schemes as a full sub-category. Iterating this
constructions as in \cite{s1}, one can also construct the homotopy
category of \textit{geometric} $n$\textit{-stacks} (which for $n$
big enough contains the homotopy category of general algebraic
stacks as a
full sub-category ).\\
This is precisely the construction we will imitate in defining our
category of geometric stacks over a symmetric monoidal model
category.

\bigskip
\bigskip

\begin{center} \textit{Geometric stacks over symmetric monoidal model categories} \end{center}

The previous construction of the homotopy category of algebraic
stacks is purely categorical. Indeed, it starts with the symmetric
monoidal category $(\mathbb{Z}-mod,\otimes)$, of
$\mathbb{Z}$-modules. The category $Aff$ of affine schemes is then
the opposite of the category of commutative and unital monoids in
the symmetric monoidal category $(\mathbb{Z}-mod,\otimes)$, which
is a categorical notion. Furthermore, the notion of a topology on
$Aff$ is also categorical. Our goal is to extend this categorical
construction to the case where $(\mathbb{Z}-mod,\otimes)$ is
replaced by a general \textit{symmetric monoidal model category}
$\mathcal{C}$ (in the sense of \cite[\S $4$]{ho}). Of course, we
want to keep track of the \textit{homotopical} information
contained in $\mathcal{C}$ and we will therefore require our
constructions
to be invariant when replacing $\mathcal{C}$ by a Quillen equivalent symmetric monoidal model category. \\

Let us start with a base \textit{symmetric monoidal model
category} $(\mathcal{C},\otimes)$ and try to imitate the
construction of algebraic stacks we have presented. The first step
is to find a reasonable analog of the category of commutative and
unital rings. It has been known since a long time by topologists
that the correct analog of commutative rings in a
\textit{homotopical} context is the notion of
$E_{\infty}$\textit{-algebra} (see for example
\cite{ekmm,hin,sp}). This notion is a generalization of the notion
of commutative monoid adapted to the case of symmetric monoidal
model categories. In particular it makes sense to consider the
category $E_{\infty}-Alg(\mathcal{C})$, of
$E_{\infty}$\textit{-algebras} in $\mathcal{C}$. Furthermore, it
is proved in \cite{bm,hin,sp} that $E_{\infty}-Alg(\mathcal{C})$ carries
a natural model category structure\footnote{To be very precise, at
this point one needs the weaker notion of \textit{semi-model}
category, but we will neglect this technical subtlety in this
introduction.}. Then, by analogy with the case of usual algebraic
geometry, we simply define the model category $\mathcal{C}-Aff$,
of \textit{affine stacks over} $\mathcal{C}$, to be the opposite
of the model category of $E_{\infty}$-algebras in $\mathcal{C}$.
It is reasonable to denote by $Spec\, A$ the object of
$\mathcal{C}-Aff$ corresponding to a $E_{\infty}$-algebra $A$. The
reader should note that if the model structure on $\mathcal{C}$ is
trivial (i.e. equivalences are isomorphisms), then
$\mathcal{C}-Aff$ is nothing else than the usual category of
commutative and unital monoids in $\mathcal{C}$, together with the
trivial model structure. In particular, if
$\mathcal{C}=\mathbb{Z}-mod$ (endowed with the trivial model
structure), the category $\mathcal{C}-Aff$ is the usual category
of affine schemes.

Our next step is to define an analog of the \textit{Yoneda
embedding} for $\mathcal{C}-Aff$. More generally, the problem is
to find a good analog of the Yoneda embedding for a model category
$M$. Of course, as an abstract category $M$ possesses the usual
Yoneda embedding, but this construction is not suited for our
purposes as it is not an invariant of the Quillen equivalence
class of $M$ (for example, it does not induces an embedding of the
homotopy category $Ho(M)$). To solve this problem, we define a
model category $M^{\wedge}$ which takes into account and depends
on the model structure on $M$. The underlying category of
$M^{\wedge}$ is as usual the category of simplicial presheaves
$\textrm{SPr}(M)$; however, the model category structure we
consider on $M^{\wedge}$ is such that its fibrant objects are
exactly objectwise fibrant simplicial presheaves $F : M^{op}
\longrightarrow SSet$ sending weak equivalences in $M$ to weak
equivalences of simplicial sets. Technically, $M^{\wedge}$ is
defined as the left Bousfield localization of the objectwise model
structure with respect to the equivalences in $M$ (see Def.
\ref{d1}). The construction $M \mapsto M^{\wedge}$ has then the
property of sending Quillen equivalences to Quillen equivalences
(see Prop. \ref{p1}). Furthermore, using mapping spaces in the
model category $M$, we construct a functor $\underline{h} : M
\longrightarrow M^{\wedge}$, which roughly speaking sends an
objects $x$ to the simplicial presheaf $y \mapsto Map_{M}(y,x)$,
$Map_{M}(-,-)$ denoting the mapping space. This functor can be
right derived into a fully faithful functor
$\mathbb{R}\underline{h} : Ho(M) \longrightarrow Ho(M^{\wedge})$,
which will be our \textit{"homotopical" Yoneda embedding} for the
model category $M$ (see Thm. \ref{t1}).

To go further one has to introduce a good notion of
\textit{Grothendieck topology} on $\mathcal{C}-Aff$ and an
associated notion of \textit{stack}. As in the previous step, we
approach, more generally, the problem of defining what is a
homotopy meaningful Grothendieck topology $\tau$ on a general
model category $M$ and what is the associated model category of
stacks $M^{\sim,\tau}$, in such a way that for trivial model
structures (i.e. when the weak equivalences are exactly the
isomorphisms) one finds back the usual notions. For this, we
introduce a notion of \textit{"homotopy" Grothendieck topology} on
a model category using the point of view of pre-topologies (see
Def. \ref{d4}). The idea is to give the usual \textit{data} of
$\tau$-coverings at the level of the homotopy category $Ho(M)$ and
require the usual \textit{conditions} of stability with respect to
\textit{isomorphisms} and \textit{composition} in $Ho(M)$ itself
while the requirement of stability under \textit{fibred products}
in $Ho(M)$ is replaced with the requirement of stability under
\textit{homotopy fibred products}. Therefore, the data of
coverings for a topology $\tau$ on $M$ are defined in $Ho(M)$
while the conditions these data have to satisfy are given at the
"higher level" in $M$ itself. This is completely natural from an
homotopic point of view and one obtains almost, but not exactly, a
usual Grothendieck topology on $Ho(M)$. We may call the pair
$(M,\tau)$ a \textit{model site}. A stack is then naturally
defined as an object in $Ho(M^{\wedge})$ satisfying a reasonable
descent condition with respect to $\tau$-hypercoverings (see Def.
\ref{d5}). We actually define a model category of stacks
$M^{\sim,\tau}$ as a certain left Bousfield localization of the
model category $M^{\wedge}$ with respect to a set $S_{\tau}$ of
maps in $M^{\wedge}$ determined by the topology.

Let us come back to our model category of affine stacks
$\mathcal{C}-Aff$. We suppose that we have chosen a topology
$\tau$ on $\mathcal{C}-Aff$ which is sub-canonical, in the sense
that the Yoneda embedding $\mathbb{R}\underline{h}$ factors
through $Ho(\mathcal{C}-Aff^{\sim,\tau})\hookrightarrow
Ho(\mathcal{C}-Aff^{\wedge})$. Therefore,
$\mathbb{R}\underline{h}$ induces a full embedding
$\mathbb{R}\underline{h} : Ho(\mathcal{C}-Aff) \longrightarrow
Ho(\mathcal{C}-Aff^{\sim,\tau})$, and objects in the essential
image of this functor will naturally be called \textit{affine
stacks over} $\mathcal{C}$. The definition of geometric stacks
over $\mathcal{C}$ is then straightforward. One consider
simplicial objects $X_{*} : \Delta^{op} \longrightarrow
\mathcal{C}-Aff^{\sim,\tau}$, which are Segal groupoids such that
$X_{0}$ is a disjoint union of affine stacks and with $X_{1}
\rightrightarrows X_{0}$ affine $\tau$-coverings. The homotopy
colimit $|X_{*}| \in Ho(\mathcal{C}-Aff^{\sim,\tau})$ of such a
simplicial object will be called a $1$\textit{-geometric stack
over} $\mathcal{C}$ for the topology $\tau$. Iterating this
construction as in \cite{s1}, one also defines
$n$\textit{-geometric stacks over} $\mathcal{C}$. The sub-category
of $Ho(\mathcal{C}-Aff^{\sim,\tau})$, consisting of $n$-geometric
stacks for some $n$, will be our
setting to \textit{do algebraic geometry over the symmetric monoidal model category} $\mathcal{C}$. \\

The following table offers a synthesis of our construction showing
how it parallels the classical constructions in algebraic
geometry.

\bigskip

$$\begin{array}{c|c}
\underline{\text{\textbf{Algebraic Geometry over}}\;
\mathbb{Z}\text{\textbf{-mod}}}&
\underline{\text{\textbf{Algebraic Geometry over a model category}}} \\
   \\
   \\
   \text{\textsf{Base Category}}: (\mathbb{Z}-\textsl{mod},\otimes) & \text{\textsf{Base Category}}: \mathcal{C}=(\mathcal{C},\otimes) \\
  \\
 alg=\text{\textsf{Commutative algebras in}}\;  (\mathbb{Z}-mod,\otimes) & Alg= E_{\infty}-\textsf{algebras}\;  \textsf{in} \; \mathcal{C} \\
  \\
 Aff=alg^{op}=\text{\textsf{Affine Schemes over}}\; (\mathbb{Z}-mod,\otimes)  &
 \mathcal{C}-Aff:=(Alg)^{opp}=\text{\textsf{Affine stacks over}} \; \mathcal{C} \\
  \\
 Aff^{\wedge} & \mathcal{C}-Aff^{\wedge} \\
   \\
 \text{\textsf{Yoneda embedding}}: & \text{\textsf{"Homotopy" Yoneda embedding}}: \\
  Aff\hookrightarrow Aff^{\wedge} & Ho(\mathcal{C}-Aff)\hookrightarrow
 Ho(\mathcal{C}-Aff^{\wedge}) \\
 \\
  \tau: \text{\textsf{Grothendieck topology on}} Aff &  \tau: \text{\textsf{"Homotopy" topology on}}\;  \mathcal{C}-Aff\\
\\
 \text{\textsf{Category of stacks}}: & \text{\textsf{Category of "homotopy" stacks}}:\\

  Ho(Aff^{\sim,\tau}) & Ho(\mathcal{C}-Aff^{\sim,\tau}) \\

\\
 \text{\textsf{Algebraic stacks in}}\;  Ho(Aff^{\sim,\tau}): |X_{*}| &
 \text{\textsf{Geometric stacks in}}\;  Ho(\mathcal{C}-Aff^{\sim}_{\tau}): |X_{*}| \\
\end{array}
$$

\bigskip
\bigskip
\bigskip

\begin{center} \textit{Examples and applications} \end{center}

The construction outlined above of the category of $n$-geometric
stacks over a symmetric monoidal model category $\mathcal{C}$ has
found his motivations in various questions coming from algebraic
geometry, algebraic topology an the recent rich interplay between
them . Among them, we describe below those which were the most
influential for us. The first two are investigated in this paper. \\

\begin{enumerate}

\item \textit{Extended or derived moduli problems.} Some of the moduli spaces arising in
Algebraic Geometry turns out to be non smooth (e.g. the moduli
stack of vector bundles over a variety of dimension greater than
one) and this maybe considered as a non-natural phenomenon. To
overcome this difficulty, the current general philosophy (see
\cite{ko},\cite{ka}, \cite{ck1}) teach us to consider the usual
moduli spaces considered so far as a \textit{truncation} of an
\textit{extended} or \textit{derived} moduli space. The non
smoothness would then arise from the fact that one is only
considering this truncation instead of the whole object. The usual
approach to these extended moduli spaces is through $DG$-schemes
(e.g. \cite{ka,ck1}). However, it was already noticed that the
homotopy category of $DG$-schemes might be not very well suited
for the functorial point of view on derived algebraic geometry.
Quoting \cite{ck2}, \\

\textit{"Similarly to the case of the usual algebro-geometric
moduli spaces, it would be nice to characterize $RHilb$ and
$RQuot$ in terms of the representability of some functors. This is
not easy, however, as the functors should be considered on the
derived category of $dg$-schemes (with quasi-isomorphisms
inverted) and for morphisms in this localized category there is
currently no explicit description. The issue should be probably
addressed in a wider foundational context for $dg$-schemes in our
present sense by means of gluing maps which are only
quasi-isomorphisms on pairwise intersections, satisfying cocycle
conditions only up to homotopy
on triple intersection etc."} \\

Our personal way of understanding the "issue" referred to in this
quotation is by stating that $DG$\textit{-schemes should be
interpreted as geometric stacks over the symmetric monoidal}
$\infty$\textit{-category of complexes}. As a first evidence for
this, we will produce a functor
$$\Theta : Ho(DG-Sch) \longrightarrow Ho(C(k)-Aff^{\sim,\textrm{ffqc}}),$$
where $C(k)$ is the symmetric monoidal model category of complexes
of $k$-modules (for any commutative and unital ring $k$),
$Ho(DG-Sch)$ is the homotopy category of $DG$-schemes over $k$
\footnote{Using $E_{\infty}$-algebra structures, the definition of
$DG$-schemes given in \cite{ck1} can be generalized over an
arbitrary ring $k$ (see Def. \ref{d31}).} and $\textrm{ffqc}$ is a
certain extension of the faithfully flat and quasi-compact
topology from usual $k$-algebras to $E_{\infty}$-algebras in
$C(k)$. We prove furthermore that $\Theta$ takes values in the
category of geometric stacks and we conjecture it is fully
faithful.

In a forthcoming version, we will also give an interpretation in
our setting of the notion of \textit{injective resolution of $BG$}
defined in \cite{ka}.

\item \textit{Brave New Algebraic Geometry.} Since the recent progress in stable algebraic
topology that led to a satisfying theory of spectra as a monoidal
model category (see \cite{ekmm}, \cite{hss}, \cite{ly}), it has
become clear that one is actually able to do usual commutative
algebra on commutative monoid objects in these categories, the so
called \textit{brave new rings}. It seems therefore natural to try
to embed this brave new commutative algebra in a \textit{brave new
algebraic geometry} i.e. in a kind of algebraic geometry over
(structured) spectra. This could give new insights in Elliptic
Cohomology and Topological Modular Forms, theories for which the
interplay between geometry and topology already proved rich and
powerful (see for example \cite{gh}, \cite {hop}, \cite{ahs},
\cite{str}).

As an example of application of our theory to this circle of
ideas, we will use the category of stacks over the symmetric
monoidal model category of symmetric spectra in order to define
the notion of \textit{\'{e}tale $K$-theory of an $E_{\infty}$-ring
spectrum}. We are not sure to deeply understand the issue of such
a construction but it certainly gives an answer to a question
pointed out to us by P.A. Ostv\ae r. To be more precise, we will
define an \'{e}tale topology on $Sp^{\Sigma}-Aff$, the model
category of affine stacks over the model category of symmetric
spectra. Then, sending each $E_{\infty}$-ring spectrum to its
$K$-theory space (as defined for example in \cite[\S $VI$]{ekmm}) gives rise
to a simplicial presheaf
$$\begin{array}{cccc}
K  : & Sp^{\Sigma}-Aff & \longrightarrow & SSet \\
& Spec\, A & \longrightarrow & K(A),
\end{array}$$
and therefore to an object $K \in
Sp^{\Sigma}-Aff^{\sim,\textrm{\'{e}t}}$. This object is in general
not fibrant (because it does not satisfy the descent condition for
\'{e}tale hypercoverings) and therefore we define for $Spec\, A
\in Sp^{\Sigma}-Aff$, $\mathbb{K}_{\textrm{\'{e}t}}(A):=RK(Spec\,
A)$, where $RK$ is a fibrant replacement of $K$ in the model
category $Sp^{\Sigma}-Aff^{\sim,\textrm{\'{e}t}}$. The space
$\mathbb{K}_{\textrm{\'{e}t}}(A)$ comes equipped with a natural
localization morphism $K(A) \longrightarrow
\mathbb{K}_{\textrm{\'{e}t}}(A)$.

\item \textit{Higher Tannakian duality.} In the preliminary manuscript \cite{to1},
$1$-Segal (or simplicial) Tannakian categories were introduced in
order to extend to higher homotopy groups the algebraic theory of
fundamental groups. The general idea was to replace in the usual
Tannakian formalism the base symmetric monoidal category of vector
spaces by the symmetric monoidal $\infty$-category of complexes.
Furthermore, as relative algebraic geometry has found interesting
applications in the Tannakian formalism (see \cite{del}), it
should not be surprising that algebraic geometry over the
$\infty$-category of complexes is relevant to higher Tannakian
theory. As an example of this principle, we will use our notion of
geometric stacks over the symmetric monoidal  model category of
complexes over some ring $k$ in order to define the notion of
\textit{affine $\infty$-gerbes}.

We start with the symmetric monoidal model category $C(k)$ of
complexes over $k$, together with a Grothendieck topology $\tau$
on $C(k)-Aff$ that will be assumed to be sub-canonical. In
practice, the choice of the topology $\tau$ is a very important
issue, but we will avoid going into these kind of considerations
here. We consider $Gp(C(k)-Aff^{\sim,\tau})$, the category of
group objects in the model category of stacks. For each $G \in
Gp(C(k)-Aff^{\sim,\tau})$, one can form its classifying simplicial
presheaf $BG \in Ho(C(k)-Aff^{\sim,\tau})$. In the case where the
underlying stack of $G$ is affine and the morphism $G
\longrightarrow *$ is a $\tau$-covering, the classifying stack
$BG$ is a $1$-geometric stack. Stacks of the form $BG$ for $G$
satisfying the above conditions will be called \textit{neutral
affine gerbes over} $C(k)$, or \textit{neutral affine}
$\infty$\textit{-gerbes over} $k$ (this definition depends on the
topology $\tau$). As the usual neutral Tannakian duality study
neutral affine gerbes (see \cite{sa}), neutral affine gerbes over
$C(k)$ are the basic object of study of higher Tannakian
duality. In the future, the higher Tannakian formalism will be
developed consistently as a certain kind of \textit{algebraic geometry over} $C(k)$.

\end{enumerate}

\bigskip

\begin{center} \textit{Relations with other works} \end{center}

The first work we would like to mention is K. Behrend's recent work on differential
graded schemes (see \cite{be}). We learned about it in one of his talks during fall
$2000$ at the MPI in Bonn. Though our interpretation of $DG$-schemes, as geometric stacks over
the model category of complexes, is similar to his own approach, the two works
seem totally independent and it is not clear to us
how the two approaches can be compared and to which extent they are really equivalent.

There are also some relations with several works on
$E_{\infty}$-algebras, in which already some standard geometrical
constructions were investigated, as for example the cotangent
complex, the tangent Lie algebra, the $K$-theory and Hochschild
cohomology spectrum, Andr\'{e}-Quillen cohomology etc. (see
\cite{ekmm,hin,gh}). We are quite convinced that all these
constructions can be generalized naturally to our setting of
geometric stacks over general model categories and will allow in
future to talk about the cotangent complex, the Lie algebra, the
cohomology or $K$-theory, Andr\'{e}-Quillen cohomology etc. of a
general geometric stack.

We have already mentioned at the beginning of this introduction
that our approach is a first approximation of what we think
algebraic geometry over symmetric monoidal $\infty$-categories
should be. Using the theory of Segal categories introduced by Z.
Tamsamani and C. Simpson (\cite{ta,sh}), it is possible to develop
such a theory without any considerations on model categories.
However, in order to compare construction of the category of
geometric stacks presented in this paper to a purely
$\infty$-categorical construction, one needs very strong version
of \textit{strictification results} (as for example in \cite[\S
$18$]{sh}). These results are already partially proved, and are
part of the foundational development of the theory of higher
categories. There is no doubt that the combined two approaches,
together with a comparison theorem allowing to pass from the world
of $\infty$-categories to the world of model categories, will be a
very powerful tool, allowing much more naturality and
manageability. As an example, let us mention that a first
consequence of such a unified theory would be a
$\infty$-categorical interpretation of the theory of
$E_{\infty}$-algebras as \textit{commutative monoids in symmetric
monoidal $\infty$-categories}. Such considerations appeared
already in T. Leisnter's work on up-to-homotopy monoid structures
(\cite{le}) and in the first author's preprint \cite{to1}.

As we have already stressed, there are some applications of our theory to the conjectural higher Tannakian formalism
described in \cite{to1}. In particular, the theory of affine stacks and schematic homotopy types of \cite{to2}, as
well as its application to non-abelian Hodge theory in \cite{kpt}, can be interpreted in terms
of algebraic geometry over the model category of complexes (at least in characteristic zero).

Finally, as explained in his letter \cite{ma}, Y. Manin's idea of
a "secondary quantization of algebraic geometry" seems to be part
of algebraic geometry over the symmetric monoidal model category
of motives (for example that defined in \cite{sp}) but our
ignorance of this subject does not allow us to say more. However,
using our notion of geometric stacks over a model category, we are
able to define an interesting candidate for the motive of an
algebraic stack (in the sense of Artin) as a \textit{$1$-geometric
stack over the model category of motives}. This construction,
which was suggested to us by the lecture of \cite{ma}, might be
closely related to the subject of secondary quantized algebraic
geometry and will be hopefully investigated in a future work.

\newpage

\begin{center} \textit{Organization of the paper} \end{center}

In the first Section of the paper we develop the theory of "homotopical" Grothendieck topologies over model
categories and the associated theory of stacks. We first define the Yoneda embedding
of a model category, then we introduce the notion of topology and construct the model category
of stacks.

In the second Section, we define and investigate the notion of geometric stack over a symmetric monoidal model
category. For this, we apply the theory of stacks developed in the first Section
to the model category of $E_{\infty}$-algebras
in a base symmetric monoidal model category. We define inductively the notion of $n$-geometric stack and
give a characterization by means of Segal groupoids as explained in this introduction.

Finally, in the third Section, we give two applications of the present theory of algebraic
geometry  over model categories.
We first explain how $DG$-schemes may be interpreted as geometric stacks over the model
category of complexes and we conclude by defining of
the \'{e}tale $K$-theory space of an $E_{\infty}$-ring spectra. \\

\begin{center} \textit{Acknowledgements} \end{center}

First, we would like to thank very warmly Markus Spitzweck for a very
exciting discussion we had with him in Toulouse a year ago, which
turned out to be the starting point of our work. We wish
especially to thank Carlos Simpson for precious conversations and
friendly encouragement: the debt we owe to his deep work on higher
categories and higher stacks will be clear throughout this work.
We are very thankful to Yuri Manin for many motivating questions on
the subject and in particular for his letter \cite{ma}. Thanks to
him, we were delighted to discover how much mathematics lies
behind the question \textit{"What is the motive of $B\mathbb{Z}/2$
?"}. For many comments and discussions, we also thank Kai
Behrend, Peter May, John Rognes and Paul-Arne Ostv\ae r. It was Paul-Arne who  pointed out to us the possible relevance of defining
\'{e}tale $K$-theory of ring spectra.

The second author wishes to thank the
Max Planck Institut f\"{u}r Mathematik in Bonn and the Laboratoire
J. A. Dieudonn\'e of the University of Nice for providing a
particularly stimulating atmosphere during his visits when part of
this work was conceived, written and partly tested in a seminar.
In particular, Andr\'{e} Hirschowitz's
enthusiasm was positive and contagious.

\end{section}

\newpage

\textit{Notations and conventions:} \\

Throughout all this work, $\mathbb{U}$ and $\mathbb{V}$ will be two universes, with $\mathbb{U} \in \mathbb{V}$, and
we will assume that $\mathbb{U}$ contains the set of natural integers, $\mathbb{N} \in \mathbb{U}$. We will use the
expression $\mathbb{U}$-set (resp. $\mathbb{U}$-group, $\mathbb{U}$-simplicial set, \dots) to denote
sets (resp. groups, resp. simplicial sets \dots) belonging to $\mathbb{U}$. The corresponding categories
will be denoted by $\mathbb{U}-Set$, $\mathbb{U}-Gp$, $\mathbb{U}-SSet$ \dots.
The words \textit{set} (resp.
\textit{group}, resp. \textit{simplicial set} \dots) will always refer to sets (resp. groups, resp. simplicial sets \dots)
belonging to the universe $\mathbb{V}$. The corresponding categories will simply be denoted by
$Set$, $Gp$, $SSet$ \dots.

We will make the following exceptions when referring to categories. A $\mathbb{U}$-category (resp.
a $\mathbb{V}$-category) will refer to a category $C$ such that for every pair of objects $(X,Y)$ in $C$, the
set $Hom(X,Y)$ belongs to $\mathbb{U}$ (resp. $\mathbb{V}$). By convention, all categories will be
$\mathbb{V}$-categories. We will say that a category is $\mathbb{U}$-small (resp. $\mathbb{V}$-small) if
it belongs to $\mathbb{U}$ (resp. to $\mathbb{V}$). \\

Our references for model categories are \cite{ho,hi}. For the weaker notion of semi-model category we refer to
\cite{sp}. An opposite category of a semi-model category will again be called a semi-model category. We will not
make a difference between the original notion and its dual.

For any simplicial semi-model category $M$, we will denote by $\underline{Hom}_{M}$ its simplicial $Hom$ set. It
will also be denoted simply by $\underline{Hom}$ when the reference to $M$ is clear. The derived version of these
simplicial $Hom$ will be denoted by $\mathbb{R}\underline{Hom}$ (see \cite[Thm. $4.3.2$]{ho}). The set of morphisms
in the homotopy category $Ho(M)$ will be denoted by $[-,-]_{M}$, or simply by $[-,-]$ when the
reference to $M$ is clear.

For a general semi-model category $M$, its mapping complexes will
be denoted by $Map_{M}$ (or $Map$ when the reference to $M$ is
clear), and will always be considered in the homotopy category of
simplicial sets (see \cite[$5.5.4$]{ho}, \cite[\S $18$]{hi},
\cite[$I.2$]{sp}). The homotopy fibred products in $Ho(M)$ will be
denoted by $x\times^{h}_{z}y$. In the same vein, the homotopy
cofibered products will be denoted by $x\coprod^{h}_{z}y$ (see
\cite[\S $11$]{hi}).

By the expression \textit{$\mathbb{V}$-cellular model categories} (resp. \textit{$\mathbb{V}$-combinatorial) model
categories}) we mean a model category satisfying the conditions of definitions \cite{hi} (resp. \cite{sm}),
expect that all sets have to be understood as $\mathbb{V}$-sets (in particular the ordinals appearing in the definition
belong to $\mathbb{V}$). \\

As usual, the standard simplicial category will be denoted by $\Delta$. For any simplicial object
$F \in \mathcal{C}^{\Delta^{op}}$ in a category $\mathcal{C}$, we will use the notation $F_{n}:=F([n])$. Similarly,
for any co-simplicial object $F \in \mathcal{C}^{\Delta}$, we will use the notation $F_{n}:=F([n])$.

\newpage

\begin{section}{Stacks over model categories}

In this first section, we will present a theory of stacks over
(semi-)model categories (we will be using \cite[\S 2]{sp} as a
reference for semi-model categories). For this, we will start by
defining the Yoneda embedding of a model category, whose idea
essentially goes back to some fundamental work of B. Dwyer and D.
Kan (see \cite{dk2}). We do not claim any originality in this
first paragraph, and the results stated are probably well known.
Then, we introduce the notion of a Grothendieck topology on a
model category, which as far as we know is a new notion. The
definition we give is very close to the usual one, and essentially
one only needs to replace in the usual definition isomorphisms by
equivalences and fibred products by homotopy fibred products.
Using this notion, we define homotopy hypercovers, which are a
straightforward generalization of hypercovers in Grothendieck's
sites,
and use them to define a model category of stacks over a model category endowed with a topology.\\
In this first version of the paper, we have not detailed the
standard properties of the model category of stacks, but we have
included some statements concerning homotopy sheaves and
computations of homotopy fibred products. These exactness
properties are fundamental to do elementary manipulations in the
model category of stacks.
Finally, we end the section by discussing the functoriality properties of the given constructions. \\

\textbf{Setting.} Throughout this section we will consider a semi-model category $M$,
together with a sub-semi-model category
$M_{\mathbb{U}} \subset M$. By this we mean that a morphism in $M_{\mathbb{U}}$ is an equivalence (resp. a fibration,
resp. a cofibration) if and only if it is an equivalence (resp. a fibration, resp. a cofibration) in $M$.
Moreover, we will
suppose that $M_{\mathbb{U}}$ is stable under the functorial factorization in $M$
(i.e. the functorial factorization functor
of $M$ can be chosen such that the factorization of a morphism in $M_{\mathbb{U}}$ stays in $M_{\mathbb{U}}$).

We will also assume that $M_{\mathbb{U}}$ is a $\mathbb{U}$-category, which is furthermore
a $\mathbb{V}$-small category. The typical example of such a situation the reader should keep in mind is
when $M$ is the model category of $\mathbb{V}$-simplicial sets (respectively, $\mathbb{V}$-simplicial groups,
complexes of $\mathbb{V}$-abelian groups, \dots) and $M_{\mathbb{U}}$ is the sub-category of
$\mathbb{U}$-simplicial sets (respectively, $\mathbb{U}$-simplicial groups,
complexes of $\mathbb{U}$-abelian groups, \dots). \\

We will make a systematic use of the left Bousfield localization
technique for which we refer to \cite[Ch. 3, 4]{hi}. In
particular, the following elementary result will be used very
frequently and we state it here merely for reference's
convenience.

\begin{prop}\label{fibrants}{(\cite[Prop. 3.6.1]{hi})}
Let $M$ be a model $\mathbb{V}$-category which is left
proper and $\mathbb{V}$-cellular (see  \cite[14.1]{hi}) or
$\mathbb{V}$-combinatorial (see \cite{sm}). If $S$ is a $\mathbb{V}$-set
of morphisms in $M$ and $L_{S}(M)$ denotes the left Bousfield
localization of $M$ with respect to $S$, then an object $W$ in
$L_{S}(M)$ is fibrant iff it is fibrant in $M$ and is $S$-local
i.e. for any map $f:{A}\rightarrow {B}$ in $S$, the induced map
between homotopy mapping spaces $f^{*}(W) : Map_{M}(B,W)\rightarrow
Map_{M}(A,W)$ is an equivalence in $SSet$.
\end{prop}

Note that we use here a slightly different notion of local objects
from that used in \cite[Def. $3.2.4$]{hi}: in Hirschhorn's
terminology fibrant objects in $L_{S}(M)$ are exactly what he
calls $S$-local objects in $M$.

\begin{subsection}{The Yoneda embedding for semi-model categories}

In this first paragraph we will construct the analog of the Yoneda embedding for semi-model categories. We will start
with the more general situation of a $\mathbb{V}$-small category $C$ together with a set of morphism $S$ in $C$
and define a model category $(C,S)^{\wedge}$. The model category $(C,S)^{\wedge}$ has to be
thought as a homotopy analog of the category of presheaves of sets $D^{\wedge}$ on a small category $D$, and is
constructed as a certain left Bousfield localization of the model category of simplicial presheaves on $C$.
This model category will be shown to be functorial in $(C,S)$ and
will only depend on the \textit{weak equivalence class} of $(C,S)$
(see Prop. \ref{p1}). Then, if $(C,S)$ is the semi-model category
$M_{\mathbb{U}}\subset M$ together with its equivalences, we will define
a Quillen adjunction
$$Re : M_{\mathbb{U}}^{\wedge} \longrightarrow M \qquad M_{\mathbb{U}}^{\wedge} \longleftarrow M : \underline{h}.$$
This adjunction will be shown to induce a fully faithful functor
$$\mathbb{R}\underline{h} : Ho(M_{\mathbb{U}}) \longrightarrow Ho(M_{\mathbb{U}}^{\wedge}),$$
which will be our final Yoneda embedding.  \\

Let us start with the general situation of a $\mathbb{V}$-small category $C$, together with a subset
of morphisms $S$ in $C$. Let $SPr(C)$ be the category of $\mathbb{V}$-simplicial presheaves on $C$, which
by \cite[Thm. $13.8.1$]{hi} is a model category where fibrations and equivalences are defined objectwise.
This model category is furthermore proper
and simplicial and the corresponding simplicial $Hom$ will simply be denoted by $\underline{Hom}$.
Recall that the simplicial structure on $SPr(C)$ is the data for any simplicial set $K$ and any
$F \in SPr(C)$, of simplicial presheaves $K\times F$ and $F^{K}$ defined by the following formulas
$$(K\times F)(x):=K\times F(x) \qquad (F^{K})(x):=\underline{Hom}_{SSet}(K,F(x)).$$
These formulas allow one to define the simplicial set of morphisms between two simplicial presheaf $F$ and $G$
by the formula
$$\underline{Hom}(F,G)_{n}:=Hom(\Delta^{n}\times F,G).$$
As usual, one has the natural adjunction isomorphisms
$$\underline{Hom}_{SPr(C)}(K\times F,G)\simeq \underline{Hom}_{SSet}(K,\underline{Hom}(F,G))\simeq
\underline{Hom}_{SPr(C)}(F,G^{K}),$$
for any simplicial set $K$ and simplicial presheaves $F,G \in SPr(C)$.

The model category $SPr(C)$ is also $\mathbb{V}$-cellular and $\mathbb{V}$-combinatorial,
therefore the left Bousfield localization techniques of \cite{hi} or \cite{sm} can be used to invert any
$\mathbb{V}$-set of maps.

Let $h_{-} : C \longrightarrow SPr(C)$ be the functor mapping an object $x \in C$ to the simplicial presheaf
it represents. In other words, $h_{x} : C^{op} \longrightarrow SSet$ send $y$ to the constant simplicial
set $Hom(y,x)$ of morphisms from $y$ to $x$ in $C$. We will denote by $h_{S}$ the image of the set of morphism
$S$ by the functor $h$. The reader should note that $h_{S}$ is a $\mathbb{V}$-set.

\begin{df}\label{d1}
\begin{itemize}
\item The simplicial model category $(C,S)^{\wedge}$ is the left Bousfield localization of $SPr(C)$ along the set
of maps $h_{S}$. When the set of map $S$ is clear from the context, we will simply write $C^{\wedge}$ for $(C,S)^{\wedge}$.
\item The derived simplicial $Hom$ of $(C,S)^{\wedge}$ will be denoted by
$$\mathbb{R}_{S}\underline{Hom}(-,-) : Ho((C,S)^{\wedge})^{op}\times Ho((C,S)^{\wedge}) \longrightarrow
Ho(SSet).$$
\end{itemize}
\end{df}

\bigskip

\textit{Important remark.} By definition, for $F,G \in (C,S)^{\wedge}$, one has
$\mathbb{R}_{S}\underline{Hom}(F,G)\simeq \mathbb{R}\underline{Hom}(F,RG)$, where
$RG$ is a fibrant model for $G$ in $(C,S)^{\wedge}$. This implies that in general, if $G$ is not fibrant
in $(C,S)^{\wedge}$, then the natural morphism
$$\mathbb{R}\underline{Hom}(F,G) \longrightarrow \mathbb{R}_{S}\underline{Hom}(F,G)$$
is not an isomorphism. This is why we need to mention the set $S$ in the notation $\mathbb{R}_{S}\underline{Hom}$. \\

We will call a simplicial presheaf $F \in SPr(C)$ $h_{S}$\textit{-local}, if for any
$h_{x} \longrightarrow h_{y}$ in $h_{S}$ the induced morphism
$$\mathbb{R}\underline{Hom}(h_{y},F) \longrightarrow\mathbb{R}\underline{Hom}(h_{x},F)$$
is an isomorphism. The reader is warned that this is a slightly
different notion of local object with respect to \cite[Def.
$3.2.4$]{hi}.

Now, the Yoneda lemma implies that one has $\mathbb{R}\underline{Hom}(h_{y},F)\simeq F(y)$; therefore, the
previous morphism is isomorphic, in the homotopy category of simplicial sets,
to the transition morphism $F(y) \longrightarrow F(x)$. This implies that $F$ is
an $h_{S}$-local object if and only if for any morphism $x \longrightarrow y$ in $C$ which is in $S$, the
induced morphism $F(y) \longrightarrow F(x)$ is an equivalence. Using this and
proposition \ref{fibrants} one finds immediately the following  result.

\begin{lem}\label{l2}
An object $F \in (C,S)^{\wedge}$ is fibrant if and only if it satisfies the following two conditions:
\begin{enumerate}
\item For every object $x \in C$, the simplicial set $F(x)$ is fibrant (i.e. $F$ is fibrant as an object
in $SPr(C)$);

\item For any morphism $x \longrightarrow y$ in $S$, the induced morphism
$F(y) \longrightarrow F(x)$ is an equivalence.

\end{enumerate}
\end{lem}

\textit{Proof:} It is a straightforward application of proposition \ref{fibrants}. \hfill $\Box$ \\

The previous lemma implies that the homotopy category $Ho((C,S)^{\wedge})$ can be naturally identified to the the full
sub-category of $Ho(SPr(C))$ consisting of simplicial presheaves $F : C^{op} \longrightarrow SSet$ sending morphisms of $S$ to
equivalences in $SSet$. Furthermore, the fibrant replacement functor in $(C,S)^{\wedge}$ induces a functor
$r : Ho(SPr(C)) \longrightarrow Ho((C,S)^{\wedge})$, which is a retraction of the natural inclusion
$Ho((C,S)^{\wedge}) \hookrightarrow Ho(SPr(C))$. \\

Let $(C,S)$ and $(D,T)$ be two $\mathbb{V}$-small categories with distinguished subsets of morphisms and $f : C \longrightarrow D$
a functor which sends $S$ into $T$. The functor induces a direct image functor on the categories of simplicial
presheaves
$$f_{*} : SPr(D) \longrightarrow SPr(C),$$
defined by $f_{*}(F)(x):=F(f(x))$, for $F \in SPr(D)$ and $x \in C$. This functor has a left adjoint
$$f^{*} : SPr(C) \longrightarrow SPr(D),$$
characterized by the property that $f^{*}(h_{x})\simeq h_{f(x)}$, for any $x \in C$.

\begin{lem}\label{l3}
The adjunction
$$f^{*} : (C,S)^{\wedge} \longrightarrow (D,T)^{\wedge} \qquad (C,S)^{\wedge} \longleftarrow (D,T)^{\wedge} : f_{*}$$
is a Quillen adjunction.
\end{lem}

\textit{Proof:} It is clear that the functor $f_{*}$ preserves levelwise equivalences and fibrations; therefore,
the adjunction $(f^{*},f_{*})$ is a Quillen adjunction between the model categories
$SPr(C)$ and $SPr(D)$. Using the general properties of left Bousfield localization of model categories (see \cite[Ch. 3, 4]{hi}),
it is then enough to prove that the functor
$f_{*}$ sends fibrant objects in $(D,T)^{\wedge}$ to fibrant objects in $(C,S)^{\wedge}$. But this is clear by
lemma \ref{l2} and the fact that $f(S)\subset T$. \hfill $\Box$ \\

\begin{df}\label{d2}
A functor $f : (C,S) \longrightarrow (D,T)$ between two $\mathbb{V}$-small categories with subsets of morphisms
is a \emph{weak equivalence} if it satisfies the following three conditions:
\begin{itemize}
\item $f(S) \subset T$;

\item There exists a functor $g : D \longrightarrow C$ with $g(T) \subset S$ and natural transformations
$$fg \Leftarrow A \Rightarrow Id \qquad gf \Leftarrow B \Rightarrow Id,$$
with $A$ (resp. $B$) an endofunctor of $D$ (resp. of $C$);

\item For any object $y \in D$ (resp. $x \in C$), the induced morphisms
$$fg(y) \longleftarrow A(y) \longrightarrow y \qquad (resp. \quad gf(x) \longleftarrow B(x) \longrightarrow x)$$
are in $T$ (resp. in $S$).
\end{itemize}
\end{df}

\begin{prop}\label{p1}
Let $f : (C,S) \longrightarrow (D,T)$ be a weak equivalence between $\mathbb{V}$-small categories
with subsets of morphisms and $g : D \longrightarrow C$ be a functor
like in defintion \ref{d3}. Then, the two Quillen adjunctions
$$f^{*} : (C,S)^{\wedge} \longrightarrow (D,T)^{\wedge} \qquad (C,S)^{\wedge} \longleftarrow (D,T)^{\wedge} : f_{*},$$
$$g^{*} : (D,T)^{\wedge} \longrightarrow (C,S)^{\wedge} \qquad (D,T)^{\wedge} \longleftarrow (C,S)^{\wedge} : g_{*},$$
are Quillen equivalences.
\end{prop}

\textit{Proof:} We will prove that the induced functor
$$\mathbb{L}f^{*} : Ho((C,S)^{\wedge}) \longrightarrow Ho((D,T)^{\wedge})$$
is an equivalence of categories, with quasi-inverse $\mathbb{L}g^{*}$.
This will be enough to show that $(f^{*},f_{*})$ and $(g^{*},g_{*})$ are Quillen equivalences.

If $F \in Ho((C,S)^{\wedge})$, let us prove that the natural morphisms
$$\mathbb{L}g^{*}\mathbb{L}f^{*}(F) \longleftarrow \mathbb{L}B^{*}(F) \longrightarrow F$$
are isomorphisms in $Ho((C,S)^{\wedge})$.
One can find a simplicial object $L_{*}$ of $(C,S)^{\wedge}$, such
that for any $[n] \in \Delta$, $L_{n}$ is isomorphic to a coproduct of simplicial presheaves of the form
$h_{x}$ with $x \in C$, together with an isomorphism in $Ho((C,S)^{\wedge})$
$$F\simeq hocolim_{[n] \in \Delta}L_{n}.$$
Then, since $\mathbb{L}g^{*}$, $\mathbb{L}f^{*}$ and $\mathbb{L}B^{*}$ commute with homotopy colimits
(because $g^{*}$, $f^{*}$ and $B^{*}$ are
left Quillen functors), one is reduced to the case where $F=h_{x}$, for some $x \in C$. But, as objects of the form
$h_{z}$ are always cofibrant, one has
$$\mathbb{L}g^{*}\mathbb{L}f^{*}(h_{x})\simeq h_{gf(x)} \qquad \mathbb{L}B^{*}(h_{x})\simeq h_{B(x)},$$
and it remains to show that the natural morphism $h_{gf(x)} \longleftarrow h_{B(x)} \longrightarrow h_{x}$ is an isomorphism in
$Ho((C,S)^{\wedge})$.
But this is true by definition of the model structure on $(C,S)^{\wedge}$ and by the fact that the morphisms
$$gf(x) \longleftarrow B(x) \longrightarrow x$$
belong to $S$.

In the same way, we prove that for any $F \in Ho((D,T)^{\wedge})$, the morphisms
$$\mathbb{L}f^{*}\mathbb{L}g^{*}(F) \longleftarrow \mathbb{L}A^{*}(F) \longrightarrow F$$
are isomorphisms in $Ho((C,S)^{\wedge})$. \hfill $\Box$ \\

\textit{Remark.} In their paper \cite{dk1}, Dwyer and Kan proved that the model category $(C,S)^{\wedge}$ is an
invariant up to Quillen equivalences of the simplicial localization
category $L(C,S)$. Proposition \ref{p1} is only a special case of this result. \\

We now come back to the basic setting of this section i.e. to our semi-model categories $M_{\mathbb{U}} \subset M$. The set of equivalences in $M_{\mathbb{U}}$
will be denoted by $\mathbf{W}_{\mathbb{U}}$.

\begin{df}\label{d3}
Let $M_{\mathbb{U}}^{c}$ (resp. $M_{\mathbb{U}}^{f}$, resp. $M_{\mathbb{U}}^{cf}$) be the sub-category of
$M_{\mathbb{U}}$ consisting of cofibrant (resps. fibrant, resp. cofibrant and fibrant) objects. We will note
$$M_{\mathbb{U}}^{\wedge}:=(M_{\mathbb{U}},\mathbf{W}_{\mathbb{U}})^{\wedge} \qquad
(M_{\mathbb{U}}^{c})^{\wedge}:=(M_{\mathbb{U}}^{c},\mathbf{W}_{\mathbb{U}}\cap M_{\mathbb{U}}^{c})^{\wedge}$$
$$(M_{\mathbb{U}}^{f})^{\wedge}:=(M_{\mathbb{U}}^{f},\mathbf{W}_{\mathbb{U}}\cap M_{\mathbb{U}}^{f})^{\wedge} \qquad
(M_{\mathbb{U}}^{cf})^{\wedge}:=(M_{\mathbb{U}}^{cf},\mathbf{W}_{\mathbb{U}}\cap M_{\mathbb{U}}^{cf})^{\wedge}.$$
These are simplicial model categories and the corresponding derived simplicial $Hom$ will be denoted by
$$\mathbb{R}_{w}\underline{Hom}(-,-) \qquad \mathbb{R}_{w,c}\underline{Hom}(-,-) \qquad
\mathbb{R}_{w,f}\underline{Hom}(-,-) \qquad \mathbb{R}_{w,cf}\underline{Hom}(-,-).$$
\end{df}

\medskip

\begin{lem}\label{l4}
The natural inclusions
$$i_{c} : M_{\mathbb{U}}^{c} \subset M_{\mathbb{U}} \qquad i_{f} : M_{\mathbb{U}}^{f} \subset M_{\mathbb{U}}  \qquad
i_{cf} : M_{\mathbb{U}}^{cf} \subset M_{\mathbb{U}},$$
induce equivalences of categories
$$\mathbb{R}(i_{c})_{*} : Ho(M_{\mathbb{U}}^{\wedge}) \simeq Ho((M_{\mathbb{U}}^{c})^{\wedge}) \qquad
\mathbb{R}(i_{f})_{*} : Ho(M_{\mathbb{U}}^{\wedge}) \simeq Ho((M_{\mathbb{U}}^{f})^{\wedge}) $$
$$\mathbb{R}(i_{cf})_{*} : Ho(M_{\mathbb{U}}^{\wedge}) \simeq Ho((M_{\mathbb{U}}^{cf})^{\wedge}).$$
These equivalences are furthermore compatible with derived simplicial $Hom$, in the sense that there
exist natural isomorphisms
$$\mathbb{R}_{w,c}\underline{Hom}(\mathbb{R}(i_{c})_{*}(-),\mathbb{R}(i_{c})_{*}(-))\simeq
\mathbb{R}_{w}\underline{Hom}(-,-) $$
$$\mathbb{R}_{w,f}\underline{Hom}(\mathbb{R}(i_{f})_{*}(-),\mathbb{R}(i_{f})_{*}(-))\simeq
\mathbb{R}_{w}\underline{Hom}(-,-)$$
$$\mathbb{R}_{w,cf}\underline{Hom}(\mathbb{R}(i_{cf})_{*}(-),\mathbb{R}(i_{cf})_{*}(-))\simeq
\mathbb{R}_{w}\underline{Hom}(-,-).$$
\end{lem}

\textit{Proof:} It is an application of proposition \ref{p1}. Let us prove for example that
$\mathbb{R}(i_{c})_{*}$ is an equivalence. For this, let $Q : M_{\mathbb{U}} \longrightarrow M_{\mathbb{U}}^{c}$ be
a cofibrant replacement functor (see \cite[p. $5$]{ho}). By definition, there exist natural transformations
$$Qi_{c} \Rightarrow Id \qquad i_{c}Q \Rightarrow Id,$$
such that for any $x \in M_{\mathbb{U}}$ the induced morphisms $Qi_{c}(x) \longrightarrow x$, $i_{c}Q(x) \longrightarrow x$
are equivalences in $M_{\mathbb{U}}$. Proposition \ref{p1} then implies that the derived functor
$\mathbb{R}(i_{c})_{*}$ is an equivalence, which preserves derived simplicial $Hom$ (as any Quillen
equivalence does).

The same proof (applied to the opposite category $M_{\mathbb{U}}^{op}$) implies that $\mathbb{R}(i_{f})_{*}$ is an
equivalence. Finally, to prove that $\mathbb{R}(i_{cf})_{*}$ is an equivalence one applies proposition \ref{p1}
first to a cofibrant replacement functor $Q : M_{\mathbb{U}} \longrightarrow M_{\mathbb{U}}^{c}$ and then
to the restriction of a fibrant replacement functor $R : M_{\mathbb{U}}^{c} \longrightarrow M_{\mathbb{U}}^{cf}$. \hfill $\Box$ \\

Lemma \ref{l4} is useful to establish functorial properties of the homotopy category
$Ho(M_{\mathbb{U}}^{\wedge})$. Indeed, if $f : M_{\mathbb{U}} \longrightarrow N_{\mathbb{U}}$ is a functor
whose restriction to $M_{\mathbb{U}}^{cf}$
preserves equivalences, then $f$ induces well defined functors
$$\mathbb{R}f_{*} : Ho(N_\mathbb{U}^{\wedge}) \longrightarrow
Ho((M_\mathbb{U}^{cf})^{\wedge})\simeq Ho(M_\mathbb{U}^{\wedge}),$$
$$\mathbb{L}f^{*} : Ho(M_\mathbb{U}^{\wedge})\simeq Ho((M_\mathbb{U}^{cf})^{\wedge}) \longrightarrow
Ho(N_\mathbb{U}^{\wedge}).$$
The functor $\mathbb{R}f_{*}$ is clearly right adjoint to $\mathbb{L}f^{*}$.

For example, let $f : M_{\mathbb{U}} \longrightarrow N_{\mathbb{U}}$ be a right Quillen functor. Then, the restriction
of $f$ on $M_{\mathbb{U}}^{cf}$ preserves equivalences and therefore induces a well defined functor
$$\mathbb{R}f_{*} : Ho(N_{\mathbb{U}}^{\wedge}) \longrightarrow Ho(M_{\mathbb{U}}^{\wedge}).$$
The same argument applies to a left Quillen functor $g : M_{\mathbb{U}} \longrightarrow N_{\mathbb{U}}$,
which by restriction to the subcategory of cofibrant and fibrant objects induces a well defined functor
$$\mathbb{R}g^{*} : Ho(N_{\mathbb{U}}^{\wedge}) \longrightarrow Ho(M_{\mathbb{U}}^{\wedge}).$$

\medskip

\begin{df}\label{d3'}
Let $M_{\mathbb{U}}$ and $N_{\mathbb{U}}$ be two $\mathbb{V}$-small semi-model categories and
$f : M_{\mathbb{U}} \longrightarrow N_{\mathbb{U}}$ be a functor whose restriction
to $M_{\mathbb{U}}^{cf}$ preserves equivalences. The previously defined functor
$$\mathbb{R}f_{*} : Ho(N_\mathbb{U}^{\wedge}) \longrightarrow Ho(M_\mathbb{U}^{\wedge})$$
will be called the \emph{inverse image} functor. Its left adjoint
$$\mathbb{L}f^{*} : Ho(M_\mathbb{U}^{\wedge}) \longrightarrow Ho(N_\mathbb{U}^{\wedge})$$
will be called the \emph{direct image} functor.
\end{df}

\textit{Remark.} The reader
should be warned that the direct and inverse image functors' construction is not functorial in $f$. In other words, if one does not
add some hypotheses on the functor $f$, then in general $\mathbb{R}f_{*}\circ \mathbb{R}g_{*}$ is not
isomorphic to $\mathbb{R}(g\circ f)_{*}$. However, one has the following easy proposition, which
ensures in many cases the functoriality of the previous construction.

\begin{prop}\label{p1'}
\begin{enumerate}
\item
Let $M_{\mathbb{U}}$, $N_{\mathbb{U}}$ and $P_{\mathbb{U}}$ be $\mathbb{V}$-small semi-model categories and
$$\xymatrix{M_{\mathbb{U}} \ar[r]^-{f} & N_{\mathbb{U}} \ar[r]^-{g} & P_{\mathbb{U}}}$$
be two functors preserving fibrant objects and equivalences between them. Then, there exist natural
isomorphisms
$$\mathbb{R}(g\circ f)_{*}\simeq \mathbb{R}f_{*}\circ \mathbb{R}g_{*} :
Ho((P_{\mathbb{U}})^{\wedge}) \longrightarrow Ho((M_{\mathbb{U}})^{\wedge}),$$
$$\mathbb{L}(g\circ f)^{*}\simeq \mathbb{L}g^{*}\circ \mathbb{L}f^{*} :
Ho((M_{\mathbb{U}})^{\wedge}) \longrightarrow Ho((P_{\mathbb{U}})^{\wedge}).$$
These isomorphisms are furthermore associative and unital in the arguments $f$ and $g$.

\item Let $M_{\mathbb{U}}$, $N_{\mathbb{U}}$ and $P_{\mathbb{U}}$ be $\mathbb{V}$-small semi-model categories and
$$\xymatrix{M_{\mathbb{U}} \ar[r]^-{f} & N_{\mathbb{U}} \ar[r]^-{g} & P_{\mathbb{U}}}$$
be two functors preserving cofibrant objects and equivalences between them. Then, there exist natural
isomorphisms
$$\mathbb{R}(g\circ f)_{*}\simeq \mathbb{R}f_{*}\circ \mathbb{R}g_{*} :
Ho((P_{\mathbb{U}})^{\wedge}) \longrightarrow Ho((M_{\mathbb{U}})^{\wedge}),$$
$$\mathbb{L}(g\circ f)^{*}\simeq \mathbb{L}g^{*}\circ \mathbb{L}f^{*} :
Ho((M_{\mathbb{U}})^{\wedge}) \longrightarrow Ho((P_{\mathbb{U}})^{\wedge}).$$
These isomorphisms are furthermore associative and unital in the arguments $f$ and $g$.
\end{enumerate}
\end{prop}

\textit{Proof:} The proof is the same as that of the usual property of
composition for derived Quillen functors (see \cite[Thm. $1.3.7$]{ho}), and is left to the reader. \hfill $\Box$ \\

Examples of functors as in the previous proposition are given by right or left Quillen functors. Therefore,
given
$$\xymatrix{M_{\mathbb{U}} \ar[r]^-{f} & N_{\mathbb{U}} \ar[r]^-{g} & P_{\mathbb{U}}}$$
a pair of right Quillen functors, one has
$$\mathbb{R}(g\circ f)_{*}\simeq \mathbb{R}f_{*}\circ \mathbb{R}g_{*}
\qquad \mathbb{L}(g\circ f)^{*}\simeq \mathbb{L}g^{*}\circ \mathbb{L}f^{*}.$$
In the same way, if $f$ and $g$ are left Quillen functors, one has
$$\mathbb{R}(g\circ f)_{*}\simeq \mathbb{R}f_{*}\circ \mathbb{R}g_{*} \qquad
\mathbb{L}(g\circ f)^{*}\simeq \mathbb{L}g^{*}\circ \mathbb{L}f^{*}.$$

\begin{prop}\label{p2}
If $f : M_{\mathbb{U}} \longrightarrow N_{\mathbb{U}}$ is a (right or left) Quillen equivalence between
$\mathbb{V}$-small semi-model
categories, then the induced functors
$$\mathbb{L}f^{*} : Ho(M_{\mathbb{U}}^{\wedge}) \longrightarrow Ho(N_{\mathbb{U}}^{\wedge}) \qquad
Ho(M_{\mathbb{U}}^{\wedge}) \longleftarrow Ho(N_{\mathbb{U}}^{\wedge}) : \mathbb{R}f_{*},$$
are equivalences, quasi-inverse of each others.
\end{prop}

\textit{Proof:} Let us prove the proposition in the case where $f$ is a right Quillen functor. The case
where $f$ is left Quillen is proved similarly.

Let $g : N_{\mathbb{U}} \longrightarrow M_{\mathbb{U}}$ be the left adjoint to
$f$; let us show that $\mathbb{L}g^{*}$ is quasi-inverse to $\mathbb{L}f^{*}$. For this, let us consider
the following two functors
$$f : M_{\mathbb{U}}^{f} \longrightarrow N_{\mathbb{U}} \qquad
RgQ : \xymatrix{ N_{\mathbb{U}} \ar[r]^{Q} & N_{\mathbb{U}}^{c} \ar[r]^{g} & M_{\mathbb{U}} \ar[r]^{R} &
M_{\mathbb{U}}^{f},}$$
where $Q$ is a cofibrant replacement functor and $R$ is a fibrant replacement functor.
The natural transformations $Id \Rightarrow R$, $Q \Rightarrow Id$ and
$gf \Rightarrow Id$, induce natural transformations
$$RgQf \Leftarrow gQf \Rightarrow gf \Rightarrow Id.$$
By hypothesis, for any $x \in M_{\mathbb{U}}^{f}$, the induced morphisms
$$RgQf(x) \longleftarrow gQf(x) \longrightarrow x$$
are equivalences. The same proof as in proposition \ref{p1} shows that
$\mathbb{L}g^{*}\mathbb{L}f^{*}$ is isomorphic to the identity. The dual argument then shows that
$\mathbb{L}f^{*}\mathbb{L}g^{*}$ is isomorphic to the identity. \hfill $\Box$ \\

To finish this paragraph, we will define a Quillen adjunction
$$Re : M_{\mathbb{U}}^{\wedge} \longrightarrow M \qquad M_{\mathbb{U}}^{\wedge} \longleftarrow M : \underline{h}$$
and show that the functor $\underline{h}$ induces a fully faithful embedding on the level of homotopy categories
(see Thm. \ref{t1}). \\

From now on, let $(\Gamma : M_{\mathbb{U}} \longrightarrow M_{\mathbb{U}}^{\Delta},i)$ be a fixed \textit{cofibrant resolution
functor} (see \cite[$17.1.3$]{hi}). This means that for any object $x \in M_{\mathbb{U}}$, $\Gamma(x)$
is a co-simplicial object
in $M_{\mathbb{U}}$, which is cofibrant for the Reedy model structure on $M_{\mathbb{U}}^{\Delta}$, together
with a natural weak equivalence $i(x) : \Gamma(x) \longrightarrow c(x)$, $c(x)$ being the constant co-simplicial object
in $M_{\mathbb{U}}$ at $x$. In the case that the semi-model category $M$ is simplicial, one can use the
standard cofibrant resolution functor $\Gamma(x):=\Delta^{*}\otimes x$.

At the level of model categories, the construction of the functor $\underline{h}$ will depend on the choice of
$\Gamma$, but after passing to the homotopy categories it will be shown that possibly different choices give the same
Yoneda embedding. \\

We define the functor $\underline{h}_{-} : M \longrightarrow SPr(M_{\mathbb{U}})$, by sending each $x \in M$
to the simplicial presheaf
$$\begin{array}{cccc}
\underline{h}_{x} : & M_{\mathbb{U}}^{op} & \longrightarrow & SSet \\
& y & \mapsto & Hom(\Gamma(y),x),
\end{array}$$
where $\Gamma(y)$ is the cofibrant resolution of $y$ induced by the functor $\Gamma$. To be more precise, the presheaf
of $n$-simplices of $\underline{h}_{x}$ is given by the formula
$$(\underline{h}_{x})_{n}(-):=Hom(\Gamma(-)_{n},x).$$

\begin{lem}\label{l5}
The functor $\underline{h} : M \longrightarrow SPr(M_{\mathbb{U}})$ is a right Quillen functor.
\end{lem}

\textit{Proof:} The fact that $\underline{h}$ is right Quillen is a direct verification and is proved in
detail in \cite[$9.5$]{du2}.  \hfill $\Box$ \\

\begin{cor}\label{c2}
The adjunction
$$Re : M_{\mathbb{U}}^{\wedge} \longrightarrow M \qquad M_{\mathbb{U}}^{\wedge} \longleftarrow M : \underline{h}$$
is a Quillen adjunction.
\end{cor}

\textit{Proof:} By the general properties of Bousfield localization of  model categories (see \cite[Ch. 3, 4]{hi}) and by lemma \ref{l5}
it is enough to show that the functor
$\underline{h}$ preserves fibrant objects. But, by definition of $\underline{h}$ this follows immediately from the
standard properties of mapping spaces (see \cite[\S $18$]{hi}) and lemma \ref{l2}. \hfill $\Box$ \\

\textit{Remark.} The reader should notice that if $(\Gamma', i')$ is another cofibrant resolution functor, then
the two derived functor $\mathbb{R}h_{-}$ and $\mathbb{R}h_{-}'$ obtained using respectively $\Gamma$
and $\Gamma'$ are naturally isomorphic. Therefore, our construction does not depend on the choice
of $\Gamma$ once one is passed to the homotopy category. \\

The main result of this first paragraph is the following one, that will play the role of the Yoneda embedding in our theory.

\begin{thm}\label{t1}
For any object $x \in Ho(M)$ which is isomorphic to an object in $Ho(M_{\mathbb{U}})$, the adjunction morphism
$$\mathbb{L}Re\mathbb{R}\underline{h}_{x} \longrightarrow x$$
is an isomorphism in $Ho(M)$.

Equivalently, the restriction of $\mathbb{R}\underline{h} : Ho(M) \longrightarrow Ho(M_{\mathbb{U}}^{\wedge})$
to the full subcategory of objects isomorphic to an object of $Ho(M_{\mathbb{U}})$, is fully faithful.
\end{thm}

\textit{Proof:} Let $x$ be a fibrant and cofibrant object in $M_{\mathbb{U}}$ and
$x \longrightarrow x_{*}$ a simplicial resolution of $x$ in $M_{\mathbb{U}}$ (see \cite[$17.1.2$]{hi}).
We consider the following two simplicial presheaves
$$\begin{array}{cccc}
h_{x_{*}} : & (M_{\mathbb{U}}^{c})^{op} & \longrightarrow & SSet \\
& y & \mapsto & Hom(y,x_{*}),
\end{array}$$
$$\begin{array}{cccc}
\underline{h}_{x_{*}} : & (M_{\mathbb{U}}^{c})^{op} & \longrightarrow & SSet \\
& y & \mapsto & Hom(\Gamma(y),x_{*}).
\end{array}$$
The augmentation $\Gamma(-) \longrightarrow c(-)$ and co-augmentation
$x \longrightarrow x_{*}$ induce a commutative diagram in $(M_{\mathbb{U}}^{cf})^{\wedge}$
$$\xymatrix{
h_{x} \ar[r]^-{a} \ar[d]_-{b} & \underline{h}_{x} \ar[d]^-{d} \\
h_{x_{*}}  \ar[r]^-{c} & \underline{h}_{x_{*}}. }$$
By the properties of mapping spaces (see \cite[\S $18$]{hi}), both  morphisms $c$ and $d$ are equivalences
in $SPr(M_{\mathbb{U}}^{c})$. Furthermore, the morphism $h_{x} \longrightarrow h_{x_{*}}$ is
isomorphic in $Ho(SPr(M_{\mathbb{U}}^{c}))$ to the induced morphism
$h_{x} \longrightarrow hocolim_{[n] \in \Delta}h_{x_{n}}$. As each morphism $h_{x} \longrightarrow h_{x_{n}}$ is
an equivalence in $(M_{\mathbb{U}}^{c})^{\wedge}$, this implies that $d$ is an
equivalence in $(M_{\mathbb{U}}^{c})^{\wedge}$. We deduce from this that
the natural morphism $h_{x} \longrightarrow \underline{h}_{x}$ is an equivalence in
$(M_{\mathbb{U}}^{c})^{\wedge}$.

Let us show how this implies that for any $x \in M_{\mathbb{U}}$,
the natural morphism $h_{x} \longrightarrow \mathbb{R}\underline{h}_{x}$ is an isomorphism
in $M_{\mathbb{U}}^{\wedge}$. Indeed, if $F$ is a fibrant object in $M_{\mathbb{U}}^{\wedge}$ and
$F_{c}$ is its restriction to $M_{\mathbb{U}}^{c}$, then one has
$$\mathbb{R}\underline{Hom}_{M_{\mathbb{U}}^{\wedge}}(\mathbb{R}\underline{h}_{x},F)\simeq
\mathbb{R}\underline{Hom}_{(M_{\mathbb{U}}^{c})^{\wedge}}(\underline{h}_{RQ(x)},F_{c})
\simeq \mathbb{R}\underline{Hom}_{(M_{\mathbb{U}}^{c})^{\wedge}}(h_{RQ(x)},F_{c})
\simeq F(RQ(x))$$
$$\simeq  \mathbb{R}\underline{Hom}_{M_{\mathbb{U}}^{\wedge}}(h_{RQ(x)},F)\simeq
\mathbb{R}\underline{Hom}(h_{x},F),$$
where $RQ(x)$ is a fibrant and cofibrant model for $x$ in $M_{\mathbb{U}}$. This shows,
by the Yoneda lemma for $Ho(M_{\mathbb{U}}^{\wedge})$,
that $h_{x} \longrightarrow \underline{h}_{x}$ is an equivalence in $M_{\mathbb{U}}^{\wedge}$.

Now, let $x \in M_{\mathbb{U}}^{f}$ and let us consider the natural morphism
in $Ho(M)$,
$$Re(h_{x}) \longrightarrow \mathbb{L}Re(\underline{h}_{x}) \longrightarrow x.$$
As $h_{x} \longrightarrow \underline{h}_{x}$ is an equivalence and $h_{x}$ is
cofibrant in $M_{\mathbb{U}}$, the first morphism $Re(h_{x}) \longrightarrow \mathbb{L}Re(\underline{h}_{x})$
is an isomorphism in $Ho(M)$. Therefore, to finish the proof of the theorem, it remains to show
that $Re(h_{x}) \longrightarrow x$ is an isomorphism in $Ho(M)$. But, by adjunction, for any
fibrant object $y \in M$, one has
$$[Re(h_{x}),y]_{M}\simeq [h_{x},\underline{h}_{y}]_{M_{\mathbb{U}}^{\wedge}}
\simeq \pi_{0}(\underline{Hom}(\Gamma(x),y))\simeq [x,y]_{M},$$
showing that $Re(h_{x}) \longrightarrow x$ is indeed an isomorphism in $Ho(M)$. \hfill $\Box$ \\

\bigskip

To finish this paragraph, let us notice that for any object $x \in M_{\mathbb{U}}$, the natural
morphism $i : \Gamma(-) \longrightarrow c(-)$ induces in the obvious manner a morphism
in $M^{\wedge}$, $h_{x} \longrightarrow \underline{h}_{x}$.

\begin{cor}\label{c2'}
For any object $x \in M_{\mathbb{U}}$, the natural morphism
$$h_{x} \longrightarrow \underline{h}_{x}$$
is an equivalence in the model category $M_{\mathbb{U}}^{\wedge}$.
\end{cor}

\textit{Proof:} This follows immediately from the proof of the previous theorem. \hfill $\Box$ \\

\end{subsection}

\begin{subsection}{Grothendieck topologies on semi-model categories}

In this paragraph, we present the notion of a Grothendieck topology on a semi-model category. The definition
is quite natural as it is formally obtained by replacing isomorphisms by equivalences
and fibred products by homotopy fibred products in the usual definition. \\

Recall that for any diagram $\xymatrix{x \ar[r]^{a} &  z &
\ar[l]_{b} y}$ in a semi-model category $M$, one can define a
homotopy fibred product $x\times^{h}_{z}y \in Ho(M)$ (see \cite[\S
$11$]{hi}). Explicitly, it is defined by
$$x\times^{h}_{z}y := x'\times_{y'}z',$$
where $\xymatrix{x' \ar[r]^{a'} & z' & \ar[l]_{b'} y'}$ is an equivalent diagram such that the two morphisms $a'$ and $b'$ are
fibrations and the objects $x'$, $y'$ and $z'$ are fibrant.
The object $x\times^{h}_{z}y$ only depends, up to a natural isomorphism in $Ho(M)$, on the equivalence
class of the diagram $\xymatrix{x \ar[r]^{a} & z & \ar[l]_{b} y}$. Furthermore it only depends, up to a non-natural
isomorphism, on the isomorphism class of the image of the diagram $\xymatrix{x \ar[r]^{a} & z & \ar[l]_{b} y}$ in $Ho(M)$.
In other words, for any diagram, $\xymatrix{x \ar[r]^{a} & z & \ar[l]_{b} y}$ in $Ho(M)$, the isomorphism class of the object
$x\times^{h}_{z}y \in Ho(M)$ is well defined. In the same way, the isomorphism classes of the two projections
$x\times^{h}_{z}y \longrightarrow x$ and $x\times^{h}_{z}y \longrightarrow y$ are well defined.

\begin{df}\label{d4}
A \emph{topology} $\tau$ on a $\mathbb{V}$-small semi-model
$\mathbb{U}$-category $M$ is the data for any object $x \in M$, of
a $\mathbb{V}$-set $Cov_{\tau}(x)$ of $\mathbb{U}$-small family
of objects in $Ho(M)/x$, called \emph{covering families} of $x$,
satisfying the following three conditions:

\begin{itemize}

\item \emph{(Stability)} For all $x \in M$ and any isomorphism $y \rightarrow x$ in $Ho(M)$,
the family $\{y \rightarrow x\}$ is in $Cov_{\tau}(x)$.

\item \emph{(Composition)} If $\{u_{i} \rightarrow x\}_{i \in I} \in Cov_{\tau}(x)$, and for any $i \in I$,
$\{v_{ij} \rightarrow u_{i}\}_{j \in J_{i}}$, the family $\{v_{ij} \rightarrow x\}_{i \in I, j\in J_{i}}$
is in $Cov_{\tau}(x)$.

\item \emph{(Homotopy base change)} For any $\{u_{i} \rightarrow x\}_{i \in I} \in Cov_{\tau}(x)$, and any morphism in $Ho(M)$,
$y \rightarrow x$, the family $\{u_{i}\times^{h}_{x}y \rightarrow y\}_{i \in I}$ is in $Cov_{\tau}(y)$.

\end{itemize}
A $\mathbb{V}$-small semi-model $\mathbb{U}$-category $M$
together with a topology $\tau$ will be called a ($\mathbb{V}$-small) \emph{semi-model}
($\mathbb{U}$-)\emph{site}.
\end{df}

\textit{Remark.} For any semi-model category $M$, one can form its
\textit{homotopy} $2$textit{-category} $\mathcal{D}^{\leq 2}(M)$
(see \cite{gz} and \cite[\S $2$]{sp}). This $2$-category is a fine
enough invariant of $M$ to be able to recover the homotopy fibred
products. It is then easy to check that the data of a topology on
$M$ only depends on the $2$-category $\mathcal{D}^{\leq 2}(M)$, up
to a $2$-equivalence. It seems to us that this is the reason why
the homotopy $2$-category of differential graded algebras is used
in \cite{be}. We warn the reader that however the homotopy
category of stacks we will define in Def. \ref{d6} depends on more
than just
$\mathcal{D}^{\leq 2}(M)$, as \textit{higher homotopies in $M$} enter in the definition. \\

Before going further in the study of topologies on semi-model categories, we would like to present
three examples.

\begin{itemize}

\item \textit{Trivial model structure.} Let $M$ be a $\mathbb{V}$-small
$\mathbb{U}$-category with the trivial model structure (i.e.
equivalences are isomorphisms and all morphisms are fibrations and
cofibrations). Then, $Ho(M)=M$ and the homotopy fibred products
are just fibred products. Therefore, a topology on the model
category $M$ in the sense of definition \ref{d4} is the same thing
as a usual Grothendieck topology on the category $M$.

\item \textit{Topological spaces.} Let us take as $M$ the model category of
$\mathbb{U}$-topological spaces, $Top$, and let us define
a topology $\tau$ in the following way.
A family of morphism in $Ho(Top)$, $\{X_{i} \rightarrow X\}_{i \in I}$,
$I \in \mathbb{U}$, is defined to be in $Cov_{\tau}(X)$ if
the induced map $\coprod_{i\in I}\pi_{0}(X_{i}) \longrightarrow \pi_{0}(X)$ is surjective. The reader will check
immediately that this defines a topology on $Top$ in the sense of definition \ref{d4}.

\item \textit{Negatively graded $CDGA$ (see \cite{be}).} Let $k$ be a field of characteristic zero and $M=CDGA_{k}^{op}$ be the
opposite model category of commutative and unital differential
graded $k$-algebras in negative degrees which belong to
$\mathbb{U}$ (see for example \cite{hin} for the description of
its model structure). Let $\tau_{0}$ be one of the usual topologies
defined on $k$-schemes (e.g. Zariski, Nisnevich, \'{e}tale, ffpf
or ffqc). Let us define a topology $\tau$ on $CDGA_{k}^{op}$ in the sense of Def. \ref{d4}, as follows. A family of morphisms in $Ho(CDGA_{k})$, $\{B
\rightarrow A_{i}\}_{i \in I}$, $I \in \mathbb{U}$, is defined to
be in $Cov_{\tau}(B)$ if it satisfies the two following
conditions:
\begin{enumerate}
\item The induced family of morphisms of affine $k$-schemes $$\{Spec\, H^{0}(A_{i}) \rightarrow Spec H^{0}(B)\}_{i \in I}$$
is a $\tau$-covering.

\item For any $i \in I$, one has $H^{*}(A_{i})\simeq H^{*}(B)\otimes_{H^{0}(B)}H^{0}(A_{i})$.

\end{enumerate}
The reader can check as an exercise that this actually defines a topology on the model category
$CDGA_{k}^{op}$. We will come back to this very important example in the last section of the paper.

\end{itemize}

\end{subsection}

\begin{subsection}{Homotopy hypercovers}

Using Reedy model structures (\cite[5.2]{ho}) on the category of simplicial objects in a model category, we generalize the
definition of hypercovers to the case of (semi-)model sites. \\

Let $M_{\mathbb{U}}$ be a $\mathbb{V}$-small semi-model $\mathbb{U}$-category
and let us consider $sM_{\mathbb{U}}$ the category of simplicial
objects in $M_{\mathbb{U}}$. By definition,
the category $M_{\mathbb{U}}$ has all $\mathbb{U}$-limits and
all $\mathbb{U}$-colimits so that the category $sM_{\mathbb{U}}$ is naturally enriched in
$\mathbb{U}-SSet$. Recall that for $K \in \mathbb{U}-SSet$ and
$x_{*} \in sM_{\mathbb{U}}$, one has by definition
$$\begin{array}{cccc}
K\times x_{*} : & \Delta^{op} & \longrightarrow & M_{\mathbb{U}} \\
& [n] & \mapsto & \coprod_{K_{n}}x_{n}.
\end{array}$$
For $x_{*}$ and $y_{*}$ objects in $sM_{\mathbb{U}}$, we define
$$\begin{array}{cccc}
\underline{Hom}(x_{*},y_{*}) : & \Delta^{op} & \longrightarrow & \mathbb{U}-SSet \\
& [n] & \mapsto & Hom_{sM_{\mathbb{U}}}(\Delta^{n}\otimes x_{*},y_{*}).
\end{array}$$
Finally, the exponential object $y_{*}^{K}$, for $K \in \mathbb{U}-SSet$ and
$y_{*} \in sM_{\mathbb{U}}$, is characterized by the
adjunction isomorphism $Hom(K\otimes x_{*},y_{*})\simeq Hom(x_{*},y_{*}^{K})$,
for all $x_{*} \in sM_{\mathbb{U}}$. \\

The category $sM_{\mathbb{U}}$ is endowed with its Reedy structure described in
\cite[Thm. $5.2.5$]{ho} and \cite[Prop. $2.7$]{sp},
which makes it into a $\mathbb{V}$-small
semi-model $\mathbb{U}$-category. Let us recall that equivalences in $sM_{\mathbb{U}}$ are defined
to be levelwise equivalences. Recall also, that a morphism $f : x_{*} \longrightarrow y_{*}$ is defined to be a fibration
if for all $n$ the morphism induced by the inclusion $\partial \Delta^{n} \hookrightarrow \Delta^{n}$,
$$x_{*}^{\Delta^{n}} \longrightarrow x_{*}^{\partial \Delta^{n}}
\times_{y_{*}^{\partial \Delta^{n}}}y_{*}^{\Delta^{n}},$$
is a fibration in $M_{\mathbb{U}}$. \\

For a $\mathbb{U}$-simplicial set $K$, the functor
$$\begin{array}{cccc}
(-)^{K} : & sM_{\mathbb{U}}  & \longrightarrow & M_{\mathbb{U}} \\
& x_{*} & \mapsto & (x_{*}^{K})_{0},
\end{array}$$
which sends a simplicial object $x_{*}$ to the $0$-th level of
the exponential object $x_{*}^{K}$, is a right Quillen functor. Its right derived functor will be
denoted by
$$(-)^{\mathbb{R} K} : Ho(sM_{\mathbb{U}}) \longrightarrow Ho(M_{\mathbb{U}}).$$
We will consider objects of $M_{\mathbb{U}}$ as constant
simplicial objects via the constant simplicial functor
$M_{\mathbb{U}} \longrightarrow sM_{\mathbb{U}}$. In particular,
for $x \in M_{\mathbb{U}}$ and $K \in \mathbb{U}-SSet$, we will
consider the object $x^{\mathbb{R} K} \in Ho(M_{\mathbb{U}})$.

\begin{df}\label{d5}
Let $x \in M_{\mathbb{U}}$ be an object in a semi-model site $(M_{\mathbb{U}},\tau)$. A \emph{homotopy}
$\tau$\emph{-hypercover} of $x$ in $M_{\mathbb{U}}$, is a simplicial object $u_{*} \in Ho(sM_{\mathbb{U}})$,
together with a morphism $u_{*} \longrightarrow x$ in $Ho(sM_{\mathbb{U}})$,
such that for any $n\geq 0$, the natural morphism
$$u_{*}^{\mathbb{R} \Delta^{n}} \longrightarrow
u_{*}^{\mathbb{R} \partial \Delta^{n}}\times^{h}_{x^{\mathbb{R} \partial \Delta^{n}}}x^{\mathbb{R} \Delta^{n}}$$
is a $\tau$-covering in $M_{\mathbb{U}}$.
\end{df}

\end{subsection}

\begin{subsection}{The model category of stacks}

In this paragraph we will use the notion of homotopy hypercover defined previously in order to construct
the model category of stacks over a model site. Our construction is based on a recent result of D. Dugger
identifying the model category of simplicial presheaves of \cite{ja} as the left Bousfield localization of
the model category of simplicial presheaves for the trivial topology by \textit{formally inverting
hypercovers} (see \cite{du}).
By definition, our model category of stacks over a model site $(M,\tau)$ will be the
left Bousfield localization of the model category $M^{\wedge}$ by formally inverting
$\tau$-hypercovers. In this first version of the paper, we will also state without proof a generalization of Dugger's theorem
by introducing the notion of homotopy sheaves in our setting. This result is fundamental to control
elementary manipulations in the model category of stacks (as for example, homotopy fibred products). \\

We come back to the basic setting of the present section, i.e. to an inclusion of semi-model categories $M_{\mathbb{U}} \subset M$,
together with the associated Yoneda embedding
$$\mathbb{R}\underline{h} : Ho(M_{\mathbb{U}}) \longrightarrow Ho(M_{\mathbb{U}}^{\wedge})$$
defined in the first paragraph. We will suppose that $M_{\mathbb{U}}$ is endowed with a topology
$\tau$ in the sense of definition \ref{d4}. \\

We define two $\mathbb{V}$-sets of morphisms in $M_{\mathbb{U}}^{\wedge}$ in the following way. For this, recall
the functor $h : M_{\mathbb{U}} \longrightarrow M_{\mathbb{U}}^{\wedge}=SPr(M_{\mathbb{U}})$, which maps an
object $x \in M_{\mathbb{U}}$  to the constant simplicial presheaf it represents.

For any
$\mathbb{U}$-set $I$ and any family of cofibrant objects $\{x_{i}\}_{i \in I} \in (M_{\mathbb{U}}^{c})^{I}$, we consider
the following natural morphism in $SPr(M_{\mathbb{U}})$
$$\coprod_{i\in I}h_{x_{i}} \longrightarrow h_{\coprod_{i \in I}x_{i}}.$$
When $I$ varies in the set of $\mathbb{U}$-sets and the $x_{i}$'s vary in the set of object in $M_{\mathbb{U}}^{c}$, we
find a $\mathbb{V}$-set of morphisms in $M_{\mathbb{U}}^{\wedge}$.
$$S_{sum}:=\{\coprod_{i\in I}h_{x_{i}} \longrightarrow h_{\coprod_{i \in I}x_{i}} \mid
I \in \mathbb{U}, \{x_{i}\}_{i \in I} \in (M_{\mathbb{U}}^{c})^{I}\}.$$

Now, for any fibrant object $x \in M_{\mathbb{U}}^{f}$, let $HHC(x)$ be the $\mathbb{V}$-set of
simplicial objects $u_{*} \in s(M/x)$, whose image in $Ho(sM)/x$ is a
homotopy $\tau$-hypercover of $x$ in $M_{\mathbb{U}}$ (see Def. \ref{d5}). For any $u_{*} \in HHC(x)$,
$[n] \mapsto h_{u_{n}}$ is a simplicial presheaf on $M_{\mathbb{U}}$ defined by the following formula:
$$\begin{array}{cccc}
h_{u_{*}} :  & M_{\mathbb{U}}^{op} & \longrightarrow & SSet \\
& y & \mapsto & ([n] \mapsto h_{u_{n}}(y)).
\end{array}$$
The augmentation $u_{*} \longrightarrow x$ gives then a morphism of simplicial presheaves $h_{u_{*}} \longrightarrow h_{x}$.
When $x$ varies in $M_{\mathbb{U}}^{f}$ and $u_{*}$ varies in $HHC(x)$,  we find a $\mathbb{V}$-set of morphisms
in $M_{\mathbb{U}}^{\wedge}$
$$S_{hhc}:=\{h_{u_{*}} \longrightarrow h_{x} \mid x \in M_{\mathbb{U}}^{f}, u_{*} \in HHC(x)\}.$$

\begin{df}\label{d6}
The \emph{simplicial model category of stacks on} $M_{\mathbb{U}}$ \emph{for the topology} $\tau$ is the left Bousfield localization
of the simplicial model category $M_{\mathbb{U}}^{\wedge}$ along the $\mathbb{V}$-set of morphisms
$S_{\tau}:=S_{sum}\cup S_{hhc}$.
It will be denoted by $M_{\mathbb{U}}^{\sim,\tau}$, or simply by $M_{\mathbb{U}}^{\sim}$ when the
topology $\tau$ is clear.

The derived simplicial $Hom$ of $M_{\mathbb{U}}^{\sim,\tau}$ will be denoted by
$$\mathbb{R}_{w,\tau}\underline{Hom}(-,-) : Ho(M_{\mathbb{U}}^{\sim,\tau})^{op} \times Ho(M_{\mathbb{U}}^{\sim,\tau})
\longrightarrow Ho(SSet).$$
\end{df}

\bigskip

The following characterization of fibrant objects in $M_{\mathbb{U}}^{\sim,\tau}$ is an immediate application of
the general criterion in Prop. \ref{fibrants}.

\begin{lem}\label{l7}
An object $F \in M_{\mathbb{U}}^{\sim,\tau}$ is fibrant if and only if it satisfies the following four
conditions:
\begin{enumerate}
\item For any $x \in M_{\mathbb{U}}$, the simplicial set $F(x)$ is fibrant;

\item For any equivalence $y \rightarrow x$ in $M$, the induced morphism $F(x) \longrightarrow F(y)$
is an equivalence of simplicial sets;

\item For any $\mathbb{U}$-set $I$ and any family of cofibrant objects $\{x_{i}\}_{i \in I}$ in $M_{\mathbb{U}}$, the
natural morphism of simplicial sets
$$F(\coprod_{i\in I}x_{i}) \longrightarrow \prod_{i \in I}F(x_{i})$$
is an equivalence;

\item For any fibrant object $x \in M_{\mathbb{U}}^{f}$ and any simplicial object $u_{*} \in s(M_{\mathbb{U}}/x)$, whose
image in $Ho(sM_{\mathbb{U}})/x$ is a homotopy $\tau$-hypercover, the natural morphism in $Ho(SSet)$
$$F(x) \longrightarrow holim_{[n] \in \Delta}F(u_{n})$$
is an isomorphism.

\end{enumerate}
\end{lem}

\textit{Proof:} It is a direct application of proposition \ref{fibrants}. \hfill $\Box$ \\

\bigskip

>From the previous lemma we immediately deduce that the homotopy category $Ho(M_{\mathbb{U}}^{\sim,\tau})$
can be identified with the full subcategory of $Ho(SPr(M_{\mathbb{U}}))$ of simplicial presheaves satisfying
conditions $(2)$, $(3)$ and $(4)$ of lemma \ref{l7}. Furthermore, the natural inclusion
$Ho(M_{\mathbb{U}}^{\sim,\tau}) \longrightarrow Ho(SPr(M_{\mathbb{U}}))$ has a left adjoint which is
a retraction. This retraction will be denoted by $a : Ho(SPr(M_{\mathbb{U}})) \longrightarrow Ho(M_{\mathbb{U}}^{\sim,\tau})$;
note that $a^{2}$ is naturally isomorphic to $a$.

\begin{df}\label{d7}
\begin{itemize}
\item
A \emph{stack} on $M_{\mathbb{U}}$ for the topology $\tau$ is an object $F \in Ho(SPr(M_{\mathbb{U}}))$
such that the natural morphism $F \longrightarrow a(F)$ is an isomorphism.
\item
For any $F \in Ho(SPr(M_{\mathbb{U}}))$, the \emph{stack associated} to $F$ is the stack $a(F)$.
\item The topology $\tau$ is \emph{sub-canonical} if for any $x \in Ho(M_{\mathbb{U}})$,
the object $\mathbb{R}\underline{h}_{x} \in Ho(M_{\mathbb{U}}^{\wedge})$
is stack.
\end{itemize}
\end{df}

\bigskip

\textit{Remarks:} \begin{itemize}
\item If $M_{\mathbb{U}}$ is endowed with the trivial model structure, then the model category
$M_{\mathbb{U}}^{\sim,\tau}$ is Quillen equivalent to the model category of simplicial presheaves defined
by J.F. Jardine in \cite{ja}. This is proved in \cite{du}. We will also state
a more general result which, for not necessarily trivial model structures,  identifies the equivalences in $M_{\mathbb{U}}^{\sim,\tau}$ as \textit{local equivalences}
(see Thm. \ref{t1'}).
Note however, that the model structure we use is not the one defined in \cite{ja} but rather its projective analog
described in \cite[$\S 5$]{sh} and \cite{bla}.

\item When the topology $\tau$ is trivial, then the model category $M_{\mathbb{U}}^{\sim,\tau}$
is equivalent to $M_{\mathbb{U}}^{\wedge}$. In particular, a stack
for the trivial topology is a simplicial presheaf $F :
M_{\mathbb{U}}^{op} \longrightarrow SSet$ which preserves
equivalences.

\item When the topology $\tau$ is sub-canonical, one obtains a fully faithful functor
$$\mathbb{R}\underline{h} : Ho(M_{\mathbb{U}}) \longrightarrow Ho(M_{\mathbb{U}}^{\sim,\tau}),$$
which embeds the homotopy theory of $M_{\mathbb{U}}$ into the homotopy theory of stacks over $M_{\mathbb{U}}$.
\end{itemize}

\bigskip

The following criterion for the topology $\tau$ to be sub-canonical can be deduced immediately
from lemma \ref{l7}.

\begin{cor}\label{l7''}
A topology $\tau$ on a semi-model category $M_{\mathbb{U}}$ is sub-canonical if and only if for
every homotopy $\tau$-hypercover $u_{*} \longrightarrow x$ in $M_{\mathbb{U}}$, the natural
morphism
$$hocolim_{[n] \in \Delta^{op}}u_{n} \longrightarrow x$$
is an isomorphism in $Ho(M_{\mathbb{U}})$.
\end{cor}

\textit{Proof:} It is a direct application of lemma \ref{l7} and of the universal property of
homotopy colimits. \hfill $\Box$ \\

\end{subsection}

\begin{subsection}{Exactness properties of the model category of stacks}

In this paragraph we will present a more classical definition of weak equivalences
in $M_{\mathbb{U}}^{\sim,\tau}$, closer to the one used in \cite{ja}. The comparison
theorem \ref{t1'} is an (easy) extension of a result of D. Dugger. Therefore, we will not give all details and the interested reader may consult \cite{ja,du} for further materials. \\

In the whole paragraph a topology $\tau$ is fixed on $M_{\mathbb{U}}$. \\

Let us start by saying a few words on the notion of sheaves of
sets in our setting of semi-model sites. As usual, any
$\mathbb{V}$-set will be considered as a constant
$\mathbb{V}$-simplicial set (note that such simplicial sets are
always fibrant in $SSet$) and therefore any functor
$M_{\mathbb{U}}^{op} \longrightarrow Set$ will be regarded as an
object in $SPr(M_{\mathbb{U}})$.

\begin{df}
\begin{enumerate}
\item
A \emph{sheaf of sets on the model site} $(M_{\mathbb{U}},\tau)$ is
a functor $F : M_{\mathbb{U}}^{op} \longrightarrow Set$ which
is stack when considered as an object in $Ho(SPr(M_{\mathbb{U}}))$.

\item
A \emph{morphism of sheaves} (of sets) $F \longrightarrow F'$ is a natural transformation.
\end{enumerate}

\medskip

\noindent The category of sheaves (of sets) on the model site
$(M_{\mathbb{U}},\tau)$ will be denoted by
$Sh(M_{\mathbb{U}},\tau)$.
\end{df}

\bigskip

By definition, a functor $F : M_{\mathbb{U}}^{op} \longrightarrow Set$ is a sheaf when it satisfies
the following three conditions:

\begin{enumerate}

\item For every equivalence $x \rightarrow y$ in $M_{\mathbb{U}}$, the induced morphism
$F(y) \longrightarrow F(x)$ is an isomorphism;

\item For every $\mathbb{U}$-small family of objects in $M_{\mathbb{U}}$, $\{x_{i}\}_{i \in I}$, $I \in \mathbb{U}$,
the induced morphism
$$F(\coprod_{i \in I}^{h}x_{i}) \longrightarrow \prod_{i \in I}F(x_{i})$$
is an isomorphism;

\item For any fibrant object $x \in M_{\mathbb{U}}^{f}$ and any simplicial object $u_{*} \in s(M_{\mathbb{U}}/x)$, whose
image in $Ho(sM_{\mathbb{U}})/x$ is a homotopy $\tau$-hypercover, the natural morphism in $Ho(SSet)$
$$F(x) \longrightarrow lim_{[n] \in \Delta}F(u_{n})\simeq Ker\left(
F(u_{0}) \rightrightarrows F(u_{1}) \right)$$
is an isomorphism.

\end{enumerate}

\bigskip

\textit{Remark.} By the universal property of the homotopy category, the category of sheaves
on $(M_{\mathbb{U}},\tau)$ is naturally equivalent to a full sub-category of the category of
presheaves of sets $Set^{Ho(M_{\mathbb{U}})^{op}}$. However, as the topology $\tau$ on the semi-model
category $M_{\mathbb{U}}$ does \textit{not} induce in general a topology on $Ho(M_{\mathbb{U}})$,
the category $Sh(M_{\mathbb{U}},\tau)$ is a priory not a category of sheaves in the usual
sense. \\

\begin{lem}
The natural functor $Sh(M_{\mathbb{U}},\tau) \longrightarrow Ho(SPr(M_{\mathbb{U}}))$ factors through
the full sub-category of stacks $Ho(M_{\mathbb{U}}^{\sim,\tau})$ and it is fully faithful.
\end{lem}

\textit{Proof:} Let $F \in Sh(M_{\mathbb{U}},\tau)$ be a sheaf and let us consider it as an object
in the model category $M_{\mathbb{U}}^{\sim,\tau}$. The sheaf conditions together with
lemma \ref{l7} show that $F$ is a fibrant object in $M_{\mathbb{U}}^{\sim,\tau}$ and in particular that its image in
$Ho(SPr(M_{\mathbb{U}}))$ is a stack.

Furthermore, as a sheaf is always a fibrant object in $M_{\mathbb{U}}^{\sim\tau}$, one checks immediately that
for two sheaves $F$ and $G$, the set of morphisms $[F,G]$ in $Ho(M_{\mathbb{U}}^{\sim,\tau})$ is
isomorphic to the set of natural transformations between $F$ and $G$. \hfill $\Box$ \\

\bigskip

Let us consider $(M_{\mathbb{U}},triv)$, the semi-model site with trivial topology. Then, the
category $Sh(M_{\mathbb{U}},triv)$ is equivalent to the category of functors
$F : M_{\mathbb{U}}^{op} \longrightarrow Set$ sending equivalences to isomorphisms. In particular, there
exists a natural fully faithful functor
$$Sh(M_{\mathbb{U}},\tau) \longrightarrow Sh(M_{\mathbb{U}},triv).$$

\begin{lem}\label{l7'}
The natural functor
$$Sh(M_{\mathbb{U}},\tau) \longrightarrow
Sh(M_{\mathbb{U}},triv)$$ has a left adjoint $a_{0} :
Sh(M_{\mathbb{U}},triv) \longrightarrow Sh(M_{\mathbb{U}},\tau)$.
Moreover, the functor $a_{0}$ is left exact (i.e. it commutes with
finite limits).
\end{lem}

\textit{Sketch of Proof:} The idea is to imitate the usual
associated sheaf functor construction, replacing fibred products
by homotopy fibred products. The proof of the left exactness of
$a_{0}$ is
very similar to the usual one. \hfill $\Box$ \\

\textit{Remark.} We have denoted $a_{0}$ the associated sheaf functor in order to make
a difference with the associated stack functor $a$. However, we will show (see Thm. \ref{t1'}) that
they coincide when applied to a sheaf, considered as an object in $Ho(M_{\mathbb{U}}^{\wedge})$. Note also that
$a_{0}$ is only defined for presheaves $M_{\mathbb{U}}^{op} \longrightarrow Set$ sending
equivalences in $M_{\mathbb{U}}$ to isomorphisms (i.e. for sheaves for the trivial
topology on $M_{\mathbb{U}}$). \\

In the next definition, $\pi_{n}(K)$, for $K \in SSet$ and $n\geq 0$ denote the set of homotopy classes of morphisms
$\Delta^{n} \longrightarrow K$ which sends $\partial \Delta^{n}$ to a point. More precisely,
$$\pi_{n}(K):=\pi_{0}\left(
\mathbb{R}\underline{Hom}(\Delta^{n},K)\times^{h}_{\mathbb{R}\underline{Hom}(\partial \Delta^{n},K)}
K_{0} \right),$$
where $\mathbb{R}\underline{Hom}(\Delta^{n},K) \longrightarrow \mathbb{R}\underline{Hom}(\partial \Delta^{n},K)$
is induced by the restriction to $\partial \Delta^{n} \subset \partial \Delta^{n}$ and
$K_{0}  \longrightarrow \mathbb{R}\underline{Hom}(\partial \Delta^{n},K)$ is adjoint to the natural
projection $K_{0}\times \partial \Delta^{n} \longrightarrow K_{0} \longrightarrow K$.

The natural projection $$\mathbb{R}\underline{Hom}(\Delta^{n},K)\times^{h}_{\mathbb{R}\underline{Hom}
(\partial \Delta^{n},K)} K_{0} \longrightarrow K_{0}$$
induces morphisms
$\pi_{n}(K) \longrightarrow K_{0}$, which make the $\pi_{n}(K)$ for $n>0$ (resp. for $n>1$)
group objects (resp. abelian group objects) over the set of $0$-simplices $K_{0}$.

\begin{df}
\begin{enumerate}
\item Let $F \in SPr(M_{\mathbb{U}})$ be a stack for the trivial topology on $M_{\mathbb{U}}$
(i.e. $F$ sends equivalences in $M_{\mathbb{U}}$ to equivalences in $SSet$).
The \emph{homotopy groups presheaves} of $F$ are defined by
$$\begin{array}{cccc}
\pi_{n}^{pr}(F) : & M_{\mathbb{U}}^{op} & \longrightarrow & Set \\
& x & \mapsto & \pi_{n}(F(x)).
\end{array}$$
They come equipped with a natural projection $\pi_{n}^{pr}(F) \longrightarrow F_{0}$.
The associated sheaves of $\pi_{n}^{pr}(F)$ are denoted by
$$\pi_{n}^{\tau}(F):=a_{0}(\pi_{n}^{pr}(F)).$$

\item Let $F, F' \in SPr(M_{\mathbb{U}})$ be two stacks for the trivial
topology on $M_{\mathbb{U}}$. A morphism $F \longrightarrow F'$ in $SPr(M_{\mathbb{U}})$ is called
a $\pi_{*}^{\tau}$-\emph{equivalence} if for all $n\geq 0$ the following square
is cartesian in $Sh(M_{\mathbb{U}},\tau)$
$$\xymatrix{
\pi_{n}^{\tau}(F) \ar[r] \ar[d] & \pi_{n}^{\tau}(F') \ar[d] \\
a_{0}(F_{0}) \ar[r] & a_{0}(F_{0}').}$$

\end{enumerate}

\end{df}

\bigskip
\bigskip

The main theorem of this paragraph is the following generalization of the main result proved in \cite{du}. Its proof
will not be given in this version of the paper.

\begin{thm}\label{t1'}
Let $F$ and $F'$ be stacks for the trivial topology on $M_{\mathbb{U}}$ (i.e. they preserve
equivalences) and $f : F \longrightarrow F'$ be a morphism in $SPr(M_{\mathbb{U}})$.
Then, $f$ is an equivalence in $M_{\mathbb{U}}^{\sim,\tau}$ if and
only if it is a $\pi_{*}^{\tau}$-equivalence.
\end{thm}

Besides its own interest, the previous theorem implies the following corollary which will be
crucial for the study of geometric stacks in the next paragraph.

\begin{cor}\label{cex}
Let
$$\xymatrix{
F \ar[r] \ar[d] & F_{1} \ar[d] \\
F_{2} \ar[r] & F_{0} }$$
be a homotopy cartesian diagram in $SPr(M_{\mathbb{U}})$. If $F_{1}$, $F_{2}$ and $F_{0}$ are stacks
for the trivial topology on $M_{\mathbb{U}}$, then
the natural morphism
$$F \longrightarrow F_{1}\times^{h}_{F_{0}}F_{2}$$
is an isomorphism in $Ho(M_{\mathbb{U}}^{\sim,\tau})$.

In other words, the associated stack functor $a$, when restricted
to the full sub-category of stacks for the trivial topology,
commutes with homotopy fibred products.
\end{cor}

\textit{Proof:} It is an application of theorem  \ref{t1'}, lemma \ref{l7'},
the long exact sequence in homotopy for a homotopy
fibred product of simplicial sets and an extended version of the five lemma. \hfill $\Box$ \\

\textit{Remark.} The functor $a : Ho(SPr(M_{\mathbb{U}}))
\longrightarrow Ho(M_{\mathbb{U}}^{\sim,\tau})$ will not commute
with homotopy fibred products in general. Indeed, suppose that
$\tau$ is the trivial topology and let $x$ be an object of
$M_{\mathbb{U}}$. The object $h_{x} \in SPr(M_{\mathbb{U}})$
represented by $x$, is such that $a(h_{x})\simeq
\underline{h}_{x}$. Therefore, if $a$ would commute with homotopy
fibred products, the natural functor $M_{\mathbb{U}}
\longrightarrow Ho(M_{\mathbb{U}})$ would send
fibred products to homotopy fibred products, which is not the case in general. \\

\end{subsection}

\begin{subsection}{Functoriality}

We finish this first section by the standard functoriality properties of the category of stacks i.e. with
direct and inverse images functors. \\

We have seen in $\S 1.1$ that any functor $f : M_{\mathbb{U}} \longrightarrow N_{\mathbb{U}}$
between $\mathbb{V}$-small semi-model categories, whose restriction to
$M_{\mathbb{U}}^{cf}$ preserves equivalences, gives rise to a pair of adjoint functors
$$\mathbb{L}f^{*} : Ho(M_{\mathbb{U}}^{\wedge}) \longrightarrow Ho(N_{\mathbb{U}}^{\wedge}) \qquad
Ho(M_{\mathbb{U}}^{\wedge}) \longleftarrow Ho(N_{\mathbb{U}}^{\wedge}) : \mathbb{R}f_{*}.$$

\begin{df}\label{d8}
Let $f : M_{\mathbb{U}} \longrightarrow N_{\mathbb{U}}$ be a functor between two $\mathbb{V}$-small
semi-model categories with topologies $\tau_{M}$ and $\tau_{N}$, respectively.
Let us suppose that the restriction of $f$ to $M_{\mathbb{U}}^{cf}$ preserves equivalences. Then,
$f$ is said to be \emph{continuous} if the inverse image functor
$$\mathbb{R}f_{*} : Ho(N_{\mathbb{U}}^{\wedge}) \longrightarrow Ho(M_{\mathbb{U}}^{\wedge})$$
preserves the categories of stacks.
\end{df}

It is immediate to check that if $f$ is a continuous functor, then
the functor
$$\mathbb{R}f_{*} : Ho(N_{\mathbb{U}}^{\sim,\tau_{N}}) \longrightarrow
Ho(M_{\mathbb{U}}^{\sim,\tau_{M}})$$
has a left adjoint
$$\mathbb{L}(f^{*})^{\sim} : Ho(M_{\mathbb{U}}^{\sim,\tau_{M}}) \longrightarrow
Ho(N_{\mathbb{U}}^{\sim,\tau_{N}}).$$
Explicitly, it is defined by the formula
$$\mathbb{L}(f^{*})^{\sim}(F):=a(\mathbb{L}f^{*}(F)),$$
for $F \in Ho(M_{\mathbb{U}}^{\sim,\tau_{M}}) \subset
Ho(M_{\mathbb{U}}^{\wedge})$,
$a$ being the associated stack functor. \\

\begin{prop}\label{p3}

Let $(M_{\mathbb{U}},\tau_{M})$, $(N_{\mathbb{U}},\tau_{N})$ and $(P_{\mathbb{U}},\tau_{P})$
be $\mathbb{V}$-small semi-model sites.
\begin{enumerate}
\item
Let
$$\xymatrix{M_{\mathbb{U}} \ar[r]^-{f} & N_{\mathbb{U}} \ar[r]^-{g} & P_{\mathbb{U}}}$$
be two continuous functors preserving fibrant objects and equivalences between them. Then, there exist natural
isomorphisms
$$\mathbb{R}(g\circ f)_{*}\simeq \mathbb{R}f_{*}\circ \mathbb{R}g_{*} :
Ho((P_{\mathbb{U}})^{\sim,\tau_{P}}) \longrightarrow Ho((M_{\mathbb{U}})^{\sim,\tau_{M}}),$$
$$\mathbb{L}((g\circ f)^{*})^{\sim}\simeq \mathbb{L}(g^{*})^{\sim}\circ \mathbb{L}(f^{*})^{\sim} :
Ho((M_{\mathbb{U}})^{\sim,\tau_{M}}) \longrightarrow Ho((P_{\mathbb{U}})^{\sim,\tau_{P}}).$$
These isomorphisms are furthermore associative and unital in the arguments $f$ and $g$.

\item Let
$$\xymatrix{M_{\mathbb{U}} \ar[r]^-{f} & N_{\mathbb{U}} \ar[r]^-{g} & P_{\mathbb{U}}}$$
be two continuous functors preserving cofibrant objects and equivalences between them. Then, there exist natural
isomorphisms
$$\mathbb{R}(g\circ f)_{*}\simeq \mathbb{R}f_{*}\circ \mathbb{R}g_{*} :
Ho((P_{\mathbb{U}})^{\sim,\tau_{P}}) \longrightarrow Ho((M_{\mathbb{U}})^{\sim,\tau_{M}}),$$
$$\mathbb{L}((g\circ f)^{*})^{\sim}\simeq \mathbb{L}(g^{*})^{\sim}\circ \mathbb{L}(f^{*})^{\sim} :
Ho((M_{\mathbb{U}})^{\sim,\tau_{M}}) \longrightarrow Ho((P_{\mathbb{U}})^{\sim,\tau_{P}}).$$
These isomorphisms are furthermore associative and unital in the arguments $f$ and $g$.
\end{enumerate}
\end{prop}

\textit{Proof:} It is immediate from proposition \ref{p1'} and the basic properties of the
functor $a$. \hfill $\Box$ \\

The following criterion gives some examples of continuous functors.

\begin{lem}\label{l8}
Let $f : M_{\mathbb{U}} \longrightarrow N_{\mathbb{U}}$ be a right Quillen functor between
two $\mathbb{V}$-small semi-model categories with topologies $\tau_{M}$ and $\tau_{N}$, respectively.
Suppose that $f$ satisfies the following two conditions:
\begin{enumerate}

\item For any $x \in Ho(M_{\mathbb{U}})$ and any
covering family $\{u_{i} \rightarrow x\}_{i \in I} \in Cov_{\tau_{M}}(x)$,
the induced family
$$\{\mathbb{R}f(u_{i}) \rightarrow \mathbb{R}f(x)\}_{i \in I}$$
is in $Cov_{\tau_{N}}(\mathbb{R}f(x))$.

\item The functor $\mathbb{R}f : Ho(M_{\mathbb{U}}) \longrightarrow Ho(N_{\mathbb{U}})$
commutes with coproducts.

\end{enumerate}

Then, the functor $f$ is continuous.
\end{lem}

\textit{Proof:} Let $F \in Ho(N_{\mathbb{U}}^{\sim,\tau_{N}})$ be
a stack and let us prove that $\mathbb{R}f^{*}$ is a stack. For
this, we use lemma \ref{l7}. The reader should notice that
conditions $(1)$ and $(2)$ are always satisfied by
$\mathbb{R}f_{*}(F)$, if they are by $F$.

Recall that by definition of the functor
$\mathbb{R}f_{*}$, for $F \in Ho(N_{\mathbb{U}}^{\sim})$ and $x \in M_{\mathbb{U}}$, one
has a natural isomorphism
$\mathbb{R}f_{*}(F)(x)\simeq F(\mathbb{R}f(x))$ in $Ho(M_{\mathbb{U}})$.

Now, by condition $(2)$ on $f$, one has
$$\mathbb{R}f_{*}(F)(\coprod^{h}_{i \in I}x_{i})\simeq \mathbb{R}\prod_{i \in I}\mathbb{R}f_{*}(F)(x_{i}),$$
for any family of objects $\{x_{i}\}_{i \in I}$ in
$Ho(N_{\mathbb{U}})$. This show that $\mathbb{R}f_{*}(F)$
satisfies $(3)$ of lemma \ref{l7}. Furthermore, as $f$ is right
Quillen it commutes with the functors $(-)^{\mathbb{R} K}$
(introduced just before  definition \ref{d5}), for any simplicial
set $K$. Using this and condition $(1)$ on $f$, one checks
immediately that $\mathbb{R}f_{*}(F)$ satisfy
condition $(4)$ of lemma \ref{l7}. \hfill $\Box$ \\

\end{subsection}

\end{section}

\begin{section}{Stacks over $E_{\infty}$-algebras}

In this second section we present the construction of the category of \textit{geometric stacks
over a base symmetric monoidal model category}, as sketched in the Introduction. For this, we will
start by recalling the homotopy theory of $E_{\infty}$-algebras and modules over them in
general symmetric monoidal model categories. The references for this part are
the foundational papers \cite{ekmm}, \cite{km}, \cite{hin} and especially \cite{sp} where the general case is studied (in particular
the case where the monoid axiom does not hold). We will then apply our theory of stacks to (semi-)model categories
of $E_{\infty}$-algebras to give a definition of geometric stacks. In the last paragraph we give
the groupoid approach to the construction of geometric stacks. \\

\textbf{Setting.} Throughout this section we will consider a left proper
$\mathbb{V}$-cofibrantly generated
symmetric monoidal model category $\mathcal{C}$. The unit of
$\mathcal{C}$ will be denoted by $\mathbf{1}$ and will always be
assumed to be a cofibrant object. We also assume $\mathcal{C}$ satisfies assumption \cite[$9.6$]{sp}, i.e.
that the domains of the generating cofibrations of $\mathcal{C}$ are cofibrant.

We will also consider $\mathcal{C}_{\mathbb{U}} \subset \mathcal{C}$ a sub-monoidal model category. By this we mean
that $\mathcal{C}_{\mathbb{U}}$ is stable under the monoidal structure  and is a sub-model category
as already explained at the beginning of section $1$. We will assume that $\mathcal{C}_{\mathbb{U}}$ is a $\mathbb{U}$-cofibrantly
generated model category which is $\mathbb{V}$-small, and that the domains and codomains of the generating
cofibrations and trivial cofibrations in $\mathcal{C}$ belong to $\mathcal{C}_{\mathbb{U}}$.

Finally, we assume that $\mathcal{C}$ is an algebra over the model category $SSet$ (of $\mathbb{V}$-simplicial sets)
or over $C(\mathbb{Z})$ (the category of complexes of $\mathbb{V}$-abelian groups). The sub-model category
$\mathcal{C}_{\mathbb{U}}$ is then assumed to be stable under external products by $\mathbb{U}$-simplicial sets
or by complexes of $\mathbb{U}$-abelian groups. \\

\begin{subsection}{Review of operads and $E_{\infty}$-algebras}

The main reference for this paragraph is \cite[$\S 1-10$]{sp}, as well as the foundational papers
\cite{ekmm,hin,km}. The reader may also consult \cite{bm}. \\

Recall first that an \textit{operad} $\mathcal{O}$ in $\mathcal{C}$ is the data,  for each integer $n \in \mathbb{N}$, of an object
$\mathcal{O}(n) \in \mathcal{C}$, together with an action of the symmetric group $\Sigma_{n}$, a unit
$\mathbf{1} \longrightarrow \mathcal{O}(1)$
and structural morphisms
$$\mathcal{O}(k)\otimes \mathcal{O}(n_{1})\otimes \dots \otimes \mathcal{O}(n_{k}) \longrightarrow
\mathcal{O}(\sum_{i}n_{i}),$$ for all integers $k\geq 1$ and
$n_{1},\dots,n_{k}$. These structural morphisms are required to
satisfy suitable associativity, commutativity and unity rules that
the reader may find in \cite[$I.1.1$]{km}. A \textit{morphism}
between two operads $\mathcal{O}$ and $\mathcal{O}'$ in
$\mathcal{C}$ is the data of morphisms $\mathcal{O}(n)
\longrightarrow \mathcal{O}'(n)$ commuting with the unit and the
structural morphisms.
With these definitions, operads in $\mathcal{C}$ form a well defined category that will be denoted by $Op(\mathcal{C})$.  \\

Following \cite{hin}, \cite{sp} and \cite{bm}, a morphism $f : \mathcal{O} \longrightarrow \mathcal{O}'$
of operads in $\mathcal{C}$ is a fibration (resp. an equivalence) if for all $n \in \mathbb{N}$, the induced
morphism $f_{n} : \mathcal{O}(n) \longrightarrow \mathcal{O}'(n)$ is a fibration (resp.
an equivalence) in $\mathcal{C}$.

\begin{thm}{(\cite[Thm. $3.2$]{sp})}\label{t2}
The category $Op(\mathcal{C})$ of operads in $\mathcal{C}$, together with the class of fibrations and equivalences
defined above, is a cofibrantly generated semi-model category.
\end{thm}

It is important to remark that $Op(\mathcal{C}_{\mathbb{U}})$ is a sub-model category of $\mathcal{C}$, which
is $\mathbb{U}$-cofibrantly generated and $\mathbb{V}$-small. \\

The fundamental operad we are interested in, is the operad $COM$, classifying commutative and unital monoids
in $\mathcal{C}$. Explicitly, it is defined by $COM(n)=\mathbf{1}$ for any $n\geq 0$, with the trivial action of
$\Sigma_{n}$.

\begin{df}{(\cite[Def. $8.1$]{sp})}\label{d9}
A \emph{unital} $E_{\infty}$-\emph{operad} in $\mathcal{C}$ is an operad $\mathcal{O} \in Op(\mathcal{C})$
satisfying the following conditions:
\begin{itemize}

\item There exists an equivalence $u : \mathcal{O} \longrightarrow COM$;

\item The induced morphism
$$u : \mathcal{O}(0) \longrightarrow COM(0)=\mathbf{1}$$
is an isomorphism;

\item For any $n\geq 0$, the object $\mathcal{O}(n)$ is cofibrant in the
semi-model category $\mathcal{C}^{\Sigma_{n}}$ of $\Sigma_{n}$-equivariant objects in $\mathcal{C}$.

\end{itemize}
\end{df}

It is important to remark that unital $E_{\infty}$-operads in $\mathcal{C}$ always exist. This follows from
\cite[Lem. $8.2$]{sp} and our assumptions on $\mathcal{C}$ which include that $\mathbf{1}$ is cofibrant
and that $\mathcal{C}$ is left proper. \\

For any operad $\mathcal{O} \in Op(\mathcal{C})$, one can define the category of $\mathcal{O}$-\textit{algebras} in $\mathcal{C}$.
By definition, an $\mathcal{O}$-algebra is the data of an object $A \in \mathcal{C}$ and structural morphisms
$\mathcal{O}(n)\otimes A^{\otimes n} \longrightarrow A$, for all $n\geq 0$. These structural morphisms
are required to satisfy certain associativity, commutativity and unit rules that the reader may find
in \cite[$I.2.1$]{km}. A morphism of $\mathcal{O}$-algebras is the data of a morphism
$A \longrightarrow A'$ in $\mathcal{C}$, commuting with the structural morphisms. These definitions
allow to define the category of algebras over a fixed operad $\mathcal{O}$ in $\mathcal{C}$, that will be
denoted by $Alg(\mathcal{O})$.

Let $f : \mathcal{O} \longrightarrow \mathcal{O}'$ be a morphism in $Op(\mathcal{C})$. Then, there exists
a natural restriction functor
$$f_{*} : Alg(\mathcal{O}') \longrightarrow Alg(\mathcal{O}).$$
This functor has a left adjoint
$$f^{*} : Alg(\mathcal{O}) \longrightarrow Alg(\mathcal{O}').$$

As for the case of operads, a morphism $f : A \longrightarrow A'$
of $\mathcal{O}$-algebras in $\mathcal{C}$, is a fibration (resp. an equivalence) if it is a fibration (resp. an equivalence) in when considered as a morphism in $\mathcal{C}$.

\begin{thm}{(\cite[Thm. $4.7$]{sp} and \cite[Cor. $6.7$]{sp})}\label{t3}
\begin{enumerate}
\item
Let $\mathcal{O} \in Op(\mathcal{C})$ be a unital $E_{\infty}$-operad in $\mathcal{C}$.
Then the category $Alg(\mathcal{O})$ of $\mathcal{O}$-algebras in $\mathcal{C}$, together with the class of fibrations
and equivalences defined above, is a cofibrantly generated semi-model category.

\item Let $f : \mathcal{O} \longrightarrow \mathcal{O}'$ be an equivalence between two
operads in $\mathcal{C}$. If for every $n\geq 0$, $\mathcal{O}(n)$ and $\mathcal{O}'(n)$ are cofibrant in
$M^{\Sigma_{n}}$,
then the induced Quillen adjunction
$$f^{*} : Alg(\mathcal{O}) \longrightarrow Alg(\mathcal{O'}) \qquad Alg(\mathcal{O}) \longleftarrow
Alg(\mathcal{O'}) : f_{*}$$
is a Quillen equivalence.
\end{enumerate}
\end{thm}

\begin{cor}\label{c3}
The semi-model category $Alg(\mathcal{O})$ of algebras over a
unital $E_{\infty}$-algebra in $\mathcal{C}$, is independent, up to a Quillen equivalence, of the choice
of $\mathcal{O} \in Op(\mathcal{C})$.
\end{cor}

\textit{Proof:} This follows from part $(2)$ of Theorem \ref{t3} and the fact that
two unitals $E_{\infty}$-operads in $\mathcal{C}$ are isomorphic in $Ho(Op(\mathcal{C}))$ (because they are
both isomorphic to $COM$). \hfill $\Box$ \\

\bigskip

The previous corollary justifies the following definition

\begin{df}\label{d10}
The semi-model category of $E_{\infty}$-\emph{algebras} in $\mathcal{C}$ is defined
to be $Alg(\mathcal{O})$, where $\mathcal{O}$ is a unital $E_{\infty}$-operad in $\mathcal{C}$. It will be denoted by $E_{\infty}-Alg(\mathcal{C})$.

The opposite semi-model category
$(E_{\infty}-Alg(\mathcal{C}))^{op}$ will be called the semi-model category of
\emph{affine stacks} over $\mathcal{C}$ and will be denoted by $\mathcal{C}-Aff$.

An $E_{\infty}$-algebra $A$ considered as an object in $\mathcal{C}-Aff$ will be
symbolically denoted by $Spec\, A$.

The same notations and terminology will be used for the category
$\mathcal{C}_{\mathbb{U}}-Aff:=(E_{\infty}-Alg(\mathcal{C}_{\mathbb{U}}))^{op}$.
\end{df}

\bigskip
\bigskip

\textit{Remarks:}\begin{itemize}
\item  By conventions, the model category $\mathcal{C}$ is an algebra over the monoidal model category
$SSet$ of simplicial sets (i.e. is a simplicial monoidal model category) or
over the category $C(\mathbb{Z})$ of complexes of abelian groups
(i.e. is a complicial monoidal model category). There are well known and famous
$E_{\infty}$-operads in $SSet$ and $C(\mathbb{Z})$, for example the singular realizations of the
\textit{little $n$-cubes
operad} and of the \textit{linear isometries operad}, as well as their
homology complexes (see \cite[$I.5$]{km}). These operads can be transported to $\mathcal{C}$ via the
unit of the algebra structure $SSet \longrightarrow \mathcal{C}$ or $C(\mathbb{Z}) \longrightarrow \mathcal{C}$
and give rise to unital $E_{\infty}$-operads in $\mathcal{C}$. This implies that in practice,
there exist natural choices for the unital $E_{\infty}$-operad in $\mathcal{C}$.

\item It is important to note that $E_{\infty}-Alg(\mathcal{C}_{\mathbb{U}})$
is a sub-model category of $E_{\infty}-Alg(\mathcal{C})$,
which is $\mathbb{U}$-cofibrantly generated and $\mathbb{V}$-small, as soon as
the $E_{\infty}$-operad has been chosen in $\mathcal{C}_{\mathbb{U}}$.
The category $Aff(\mathcal{C}_{\mathbb{U}})$, opposite to the category
$E_{\infty}-Alg(\mathcal{C}_{\mathbb{U}})$ is really the category we will
be mostly interested in.

\item If the model structure on $\mathcal{C}$ is trivial, then so is the model structure on
$E_{\infty}-Alg(\mathcal{C})$. Actually, a unital $E_{\infty}$-operad is then automatically
isomorphic to the operad $COM$. Therefore in this case, $E_{\infty}-Alg(\mathcal{C})$ is just the trivial model category
of commutative and unital monoids in $\mathcal{C}$.

\end{itemize}

Let $\mathcal{O}$ be an operad in $\mathcal{C}$ and $A$ be an $\mathcal{O}$-algebra in $\mathcal{C}$.
An $A$-\textit{module} in $\mathcal{C}$ is the data of an object $M \in \mathcal{C}$ and structural morphisms
$\mathcal{O}(n)\otimes A^{n-1}\otimes M \longrightarrow M$. These structural morphisms
are required to satisfy certain associativity, commutativity and unit rules that the reader may find
in \cite[$I.4.1$]{km}. A \textit{morphism} of $A$-modules is the data of a morphism $M \longrightarrow M'$
in $\mathcal{C}$, commuting with the structural morphisms. These definitions
allow to define the category of modules over a fixed operad algebra $A$ over a fixed operad
$\mathcal{O}$ in $\mathcal{C}$, which will be
simply denoted by $Mod(A)$.

Let $\mathcal{O} \in Op(\mathcal{C})$ be an operad in $\mathcal{C}$ and
$f : A \longrightarrow A'$ be a morphism in $Alg(\mathcal{O})$. Then, there exists
a natural restriction functor
$$f_{*} : Mod(A') \longrightarrow Mod(A).$$
This functor has a left adjoint
$$f^{*} : Mod(A) \longrightarrow Mod(A').$$

As for the case of operads and algebras, a morphism $f : M \longrightarrow M'$
of $A$-modules in $\mathcal{C}$ is a fibration (resp. an equivalence) if it is a fibration (resp. an equivalence) when considered as a morphism in $\mathcal{C}$.

\begin{thm}{(\cite[Thm. $6.1$]{sp} and \cite[Cor. $6.7$]{sp})}\label{t4}.\\
Let $\mathcal{O} \in Op(\mathcal{C})$ be a unital $E_{\infty}$-operad in $\mathcal{C}$ and
$A \in Alg(\mathcal{O})$ a cofibrant $E_{\infty}$-algebra in $\mathcal{C}$. Then
\begin{enumerate}
\item
The category $Mod(A)$ of $A$-modules in $\mathcal{C}$, together with the classes  of fibrations
and equivalences defined above, is a cofibrantly generated model category;

\item If $f : A \longrightarrow A'$ is an equivalence between two cofibrant $E_{\infty}$-algebras
in $\mathcal{C}$, the adjunction
$$f^{*} : Mod(A) \longrightarrow Mod(A') \qquad Mod(A) \longleftarrow
Mod(A') : f_{*}$$
is a Quillen equivalence.
\end{enumerate}
\end{thm}

\bigskip
\bigskip

Again, if $\mathcal{O} \in Op(\mathcal{C}_{\mathbb{U}})$ is a unital $E_{\infty}$-algebra
and $A \in Alg(\mathcal{O})$ is an $E_{\infty}$-algebra in $\mathcal{C}_{\mathbb{U}}$, then
the category of $A$-modules in $\mathcal{C}_{\mathbb{U}}$ is a
sub-model category of $Mod(A)$. It is furthermore $\mathbb{U}$-cofibrantly generated and
$\mathbb{V}$-small. \\

The compatibility between the pushforward and pullback functors above, on the model categories
of modules, is expressed through the following base change formula.

\begin{prop}{(\cite[Prop. $9.12$]{sp})}\label{p4}.
Let
$$\xymatrix{
A \ar[r]^{f} \ar[d]_{g} & B \ar[d]^-{g'} \\
A' \ar[r]_-{f'} & B' }$$
be a homotopy co-cartesian diagram of cofibrant $E_{\infty}$-algebras in $\mathcal{C}$. Then, for any
$M \in Ho(Mod(B))$, the natural base change morphism
$$\mathbb{L}(g)^{*}\mathbb{R}f_{*}(M) \longrightarrow \mathbb{R}f'_{*}\mathbb{L}(g')^{*}(M)$$
is an isomorphism in $Ho(Mod(A'))$.
\end{prop}

\end{subsection}

\begin{subsection}{Geometric stacks over a monoidal model category}

For this paragraph, recall our basic setting of this section: an
inclusion $\mathcal{C}_{\mathbb{U}} \subset  \mathcal{C}$ of
monoidal model categories,  satisfying the conditions explained at
the beginning of this section. We will assume that the semi-model
category $\mathcal{C}_{\mathbb{U}}-Aff$ of affine stacks over
$\mathcal{C}_{\mathbb{U}}$ is endowed with a topology $\tau$. This
semi-model site will be denoted by
$(\mathcal{C}_{\mathbb{U}}-Aff,\tau)$ and the corresponding model
category of stacks will simply be denoted by
$\mathcal{C}_{\mathbb{U}}-Aff^{\sim,\tau}$. For the  sake of
simplicity we assume the topology $\tau$ is sub-canonical (see
Def. \ref{d7}).

Recall from Section $1$ the existence of a Quillen adjunction
$$\textrm{Re} : \mathcal{C}_{\mathbb{U}}-Aff^{\sim,\tau} \longrightarrow \mathcal{C}-Aff
\qquad \mathcal{C}_{\mathbb{U}}-Aff^{\sim,\tau} \longleftarrow \mathcal{C}-Aff :
\texttt{Spec},$$
where we denote by $\texttt{Spec}$ the functor we have called $\underline{h}$ in Section 1, because this seems more natural when dealing with $E_{\infty}$-algebras. This Quillen adjunction, as we saw in Section 1, induces a fully faithful functor (the Yoneda embedding)
$$\mathbb{R}\underline{Spec} :  Ho(\mathcal{C}_{\mathbb{U}}-Aff)=Ho(E_{\infty}-Alg(\mathcal{C}_{\mathbb{U}}))^{op}
\longrightarrow Ho(\mathcal{C}_{\mathbb{U}}-Aff^{\sim,\tau}).$$

In order to give the definition of $n$-\textit{geometric stacks}, we will need
to add the following hypothesis on our topology $\tau$:\\

\begin{hyp}\label{hyp} Let $\{Spec\, B_{i}\longrightarrow Spec\, A\}_{i \in I}$ be a
$\mathbb{U}$-small family of morphisms in $Ho(\mathcal{C}_{\mathbb{U}}-Aff)$ and, for each $i \in I$, let
$\{Spec\, C_{j}\longrightarrow Spec\, B_{i}\}_{j\in J_{i}}$ be a
$\tau$-covering family. If the induced family in $Ho(\mathcal{C}_{\mathbb{U}}-Aff)$,
$\{Spec\, C_{j}\longrightarrow Spec\, A\}_{i\in I, j \in J_{i}}$
is a $\tau$-covering, then so is $\{Spec\, B_{i}\longrightarrow Spec\, A\}_{i \in I}$.
\end{hyp}

The definition of an $n$-geometric stack over $\mathcal{C}_{\mathbb{U}}$
is given by induction on $n$. We define simultaneously the notion of $n$-geometric
stack and the notion of $n$-covering family by induction on $n$. Note that the
both definitions depend on the topology $\tau$, and one should probably use the
expression $n_{\tau}$-geometric stacks. However, as we will not consider different
topologies at the same time, we will omit the reference to $\tau$.

\begin{df}\label{d20}
\begin{itemize}
\item The category of $0$-\emph{geometric stacks} over
$\mathcal{C}_{\mathbb{U}}$ is the essential image of the functor
$\mathbb{R}\texttt{Spec}$. It will be denoted by $0-GeSt(\mathcal{C}_{\mathbb{U}})$ and
is equivalent to $\mathcal{C}_{\mathbb{U}}-Aff$ via the Yoneda embedding. Note also that it does not depend on $\tau$.

The category $0-GeSt$ will also be called the category of \emph{affine stacks} over $\mathcal{C}_{\mathbb{U}}$.

\item A morphism $f : F \longrightarrow F'$ in $Ho(\mathcal{C}_{\mathbb{U}}-Aff^{\sim,\tau})$
is $0$-\emph{representable} if for any $0$-geometric stack $H$ and any morphism
$H \longrightarrow F'$, the homotopy pull-back
$F\times^{h}_{F'}H$ is a $0$-geometric stack (this is again independent of   $\tau$).

\item A $\mathbb{U}$-small family of morphisms
$\{f_{i} : F_{i} \longrightarrow F'\}_{i\in I}$, $I\in \mathbb{U}$,
in $Ho(\mathcal{C}_{\mathbb{U}}-Aff^{\sim,\tau})$, is a $0$-\emph{covering} if it satisfies the
two following conditions:

\begin{itemize}
\item For any $i \in I$, the morphism $f_{i}$ is $0$-representable;

\item For any $0$-geometric stack $H$, any morphism
$H \longrightarrow F'$ and any $i \in I$, the homotopy pull-back family
$\{F_{i}\times^{h}_{F'}H \longrightarrow H\}_{i \in I}$ (which is
a  $\mathbb{U}$-small family of morphisms of $0$-geometric stacks by the first condition), corresponds
to a $\tau$-covering family in $Ho(\mathcal{C}_{\mathbb{U}}-Aff)$.

\end{itemize}

\end{itemize}

Let us suppose that the full sub-category $(n-1)-GeSt(\mathcal{C}_{\mathbb{U}})
\subset Ho(\mathcal{C}_{\mathbb{U}}-Aff^{\sim,\tau})$
of $(n-1)$-geometric stacks
has been defined, as well as the notion of
a $(n-1)$-covering family in $Ho(\mathcal{C}_{\mathbb{U}}-Aff^{\sim,\tau})$.

\begin{itemize}

\item A morphism $f : F \longrightarrow F'$ in $Ho(\mathcal{C}_{\mathbb{U}}-Aff^{\sim,\tau})$ is
$(n-1)$-\emph{representable} if, for every $0$-geometric stack $H$
and any morphism $H \longrightarrow F'$, the homotopy fibred
product $F\times^{h}_{F'}H \in
Ho(\mathcal{C}_{\mathbb{U}}-Aff^{\sim,\tau})$ is an
$(n-1)$-geometric stack.

\item A stack $F \in Ho(\mathcal{C}_{\mathbb{U}}-Aff^{\sim,\tau})$ is $n$-\emph{geometric} if
it satisfies the following two conditions:
\begin{itemize}
\item The diagonal morphism $F \longrightarrow F\times F$ is
$(n-1)$-representable.
\item There exists a $\mathbb{U}$-small $(n-1)$-covering family
$\{f_{i} : F_{i} \longrightarrow F\}_{i \in I}$, $I\in \mathbb{U}$, such that
each $F_{i}$ is a $0$-geometric stack. Such a family will be
called a $(n-1)$-\emph{atlas} for $F$.

\end{itemize}

The full sub-category of $Ho(\mathcal{C}_{\mathbb{U}}-Aff^{\sim,\tau})$ consisting of
$n$-geometric stacks will be denoted by $n-GeSt(\mathcal{C}_{\mathbb{U}})$.

\item Let $F$ be a $0$-geometric stack.
A $\mathbb{U}$-small family of morphisms $\{f_{i} : F_{i} \longrightarrow F\}_{i \in I}$ in $n-GeSt(\mathcal{C}_{\mathbb{U}})$
is a \emph{special} $n$-\emph{covering}
if, for any $i\in I$, there exists a $(n-1)$-atlas $\{H_{i,j} \longrightarrow
F_{i}\}_{j\in J_{i}}$,
such that the induced family $\{H_{i,j} \longrightarrow F\}_{i\in I,j\in J_{i}}$
is a $0$-covering in $Ho(\mathcal{C}_{\mathbb{U}}-Aff^{\sim,\tau})$.

\item A $\mathbb{U}$-small family of morphisms
$\{f_{i} : F_{i} \longrightarrow F\}_{i \in I}$ in $Ho(\mathcal{C}_{\mathbb{U}}-Aff^{\sim,\tau})$ is a $n$-\emph{covering}
if it satisfies the following two conditions>

\begin{itemize}
\item For any $i \in I$, the morphism $f_{i}$ is $n$-representable;

\item For any $0$-geometric stack $H$ and any morphism
$H \longrightarrow F'$, the homotopy pull-back family
$\{F_{i}\times^{h}_{F'}H \longrightarrow H\}_{i\in I}$ (which is
a family of morphisms from $n$-geometric stacks to a $0$-geometric stack),
is a special $n$-covering family (as defined before).

\end{itemize}

\end{itemize}

A stack will be simply called geometric if it is $n$-geometric for some integer $n$.

\end{df}

\bigskip
\bigskip

\textit{Remarks:}
\begin{itemize}
\item For $F$ an $n$-geometric stack, the integer $n$ refers
to the complexity of the geometry of $F$ and not to
its homotopical complexity as the usual expression \textit{$n$-stack} refers to.
In general, the notion of $n$-geometric stack has nothing to do with the
notion of $n$-stack i.e. of $n$-truncated simplicial presheaf. These two notions
relate each others only when the model structure on $\mathcal{C}$ is
trivial and the reason is that in this case affine stacks are
$0$-truncated (i.e. are presheaves of constant simplicial sets).

\item The reader should be warned that when $\mathcal{C}_{\mathbb{U}}$
is the monoidal trivial model category of $\mathbb{U}$-abelian
groups, then our notion of $n$-geometric stacks for $n=0,1$ is \textit{not}
equivalent to the notion of algebraic spaces and algebraic stacks
as commonly used (e.g., in \cite{lm}), say with $\tau$ the
$\textrm{ffqc}$-topology. For example, a non-affine scheme is not
a $0$-geometric stack in our sense but it is a $1$-geometric stack
if it is separated. To get non-separated schemes one
needs to consider $2$-geometric stacks. In the same way, an Artin
stack with a non-affine diagonal is not a $1$-geometric stack in
our sense. It is a $2$-geometric stack if it is quasi-separated
but only a $3$-geometric stack in general.

\end{itemize}

Using a big induction argument on $n$ like it is done in \cite{s1}, one proves
the following basic proposition. We will not rewrite the argument in this
version of the paper. Note however that the proof uses in an essential way the hypothesis \ref{hyp} on our topology. This hypothesis is precisely used to prove independence of the choice of atlases.

\begin{prop}\label{d22} With the notations as above:
\begin{enumerate}
\item There are natural inclusions $n-GeSt(\mathcal{C}_{\mathbb{U}})\subset (n+1)-GeSt(\mathcal{C}_{\mathbb{U}})$;

\item The set of $n$-representable morphisms is stable by
composition and base change. Any isomorphism is $n$-representable
for any $n$;

\item The set of $n$-covering families is stable by compositions and base changes. Any
isomorphism is a $n$-covering family for any $n$;

\item If $f : F \longrightarrow F'$ is a $n$-representable morphism and
$F'$ is a $n$-geometric stack, then so is $F$;

\item
The sub-category $n-GeSt(\mathcal{C}_{\mathbb{U}}) \subset
Ho(\mathcal{C}_{\mathbb{U}}-Aff^{\sim,\tau})$ is stable under
homotopy fibred products.

\end{enumerate}
\end{prop}

\textit{Proof:} See \cite{s1}. \hfill $\Box$

\end{subsection}

\begin{subsection}{An example: Quotient stacks}

The model category of stacks
$\mathcal{C}_{\mathbb{U}}-Aff^{\sim,\tau}$ is a
$\mathbb{U}$-cofibrantly generated model category and therefore,
for any $\mathbb{U}$-small category $I$, the category of
$I$-diagrams $(\mathcal{C}_{\mathbb{U}}-Aff^{\sim,\tau})^{I}$ is a
model category with the so-called projective model structure (see
\cite[Thm. $13.8.1$]{hi}). Let us recall that fibrations and
equivalences in $(\mathcal{C}_{\mathbb{U}}-Aff^{\sim,\tau})^{I}$
are defined levelwise. We will be interested in the case
$I=\Delta^{op}$ i.e. in the category of simplicial objects in
$\mathcal{C}_{\mathbb{U}}-Aff^{\sim,\tau}$. As usual we will
denote
$s\mathcal{C}_{\mathbb{U}}-Aff^{\sim,\tau}:=(\mathcal{C}_{\mathbb{U}}-Aff^{\sim,\tau})^{\Delta^{op}}$
and, for $X_{*} \in s\mathcal{C}_{\mathbb{U}}-Aff^{\sim,\tau}$,
$X_{n}:=X([n])$.

Recall that for any integer $n>0$ and any $X_{*} \in s\mathcal{C}_{\mathbb{U}}-Aff^{\sim,\tau}$,
there exists a Segal morphism
$$S_{n} : X_{n} \longrightarrow \underbrace{X_{1}\times^{h}_{X_{0}}X_{1}
\times^{h}_{X_{0}} \dots \times^{h}_{X_{0}}X_{1}}_{n\; times}$$
induced by the morphisms in $\Delta$,
$$\alpha_{i} : [1] \longrightarrow [n] \qquad d_{o} : [0] \longrightarrow [1]\qquad d_{1} : [0] \longrightarrow [1],$$
which are defined by
$$\alpha_{i}(0)=i \qquad \alpha_{i}(1)=i+1, \qquad  0\leq i<n.$$

The following definition is a generalization of Segal's $\Delta^{op}$-spaces to the case
where the \textit{space of objects} is not contractible.

\begin{df}\label{d21}
An object $X_{*} \in s\mathcal{C}_{\mathbb{U}}-Aff^{\sim,\tau}$ is called a \emph{Segal groupoid object} if it satisfies the
following two conditions:
\begin{enumerate}
\item For every integer $n>0$, the Segal morphism
$$S_{n} : X_{n} \longrightarrow \underbrace{X_{1}\times^{h}_{X_{0}}X_{1}
\times^{h}_{X_{0}} \dots \times^{h}_{X_{0}}X_{1}}_{n\; times},$$
is an equivalence in $\mathcal{C}_{\mathbb{U}}-Aff^{\sim,\tau}$;

\item The natural morphism
$$d_{2}\times d_{1} : X_{2} \longrightarrow X_{1}\times_{d_{1},X_{0},d_{1}}^{h}X_{1}$$
is an equivalence in $\mathcal{C}_{\mathbb{U}}-Aff^{\sim,\tau}$.

\end{enumerate}
\end{df}

\bigskip

\textit{Remark.} Condition $(2)$ implies that the induced simplicial object in $Ho(\mathcal{C}_{\mathbb{U}}-Aff^{\sim,\tau})$
is a groupoid object. \\

The colimit functor $colim : s\mathcal{C}_{\mathbb{U}}-Aff^{\sim,\tau} \longrightarrow \mathcal{C}_{\mathbb{U}}-Aff^{\sim,\tau}$
is clearly a left
Quillen functor and can then be left derived to a functor at the level of homotopy categories
$$|-|:=hocolim_{\Delta^{op}} : Ho(s\mathcal{C}_{\mathbb{U}}-Aff^{\sim,\tau}) \longrightarrow Ho(\mathcal{C}_{\mathbb{U}}-Aff^{\sim,\tau}).$$

\begin{df}
If $X_{*} \in Ho(s\mathcal{C}_{\mathbb{U}}-Aff^{\sim,\tau})$ is a Segal groupoid, then $|X_{*}| \in
Ho(\mathcal{C}_{\mathbb{U}}-Aff^{\sim,\tau})$ is called
the \emph{quotient stack} of $X_{*}$.
\end{df}

The fundamental theorem of this section is the following generalization of the criterion of \cite[Prop. $4.1$]{s1}. We will only provide a sketch of its
proof in this version of the paper.

\begin{thm}\label{t5}
A stack $F \in Ho(\mathcal{C}_{\mathbb{U}}-Aff^{\sim,\tau})$ is $n$-geometric if and only if there exists a Segal groupoid
$X_{*} \in Ho(s\mathcal{C}_{\mathbb{U}}-Aff^{\sim,\tau})$ such that $F\simeq |X_{*}|$ and satisfying the following
two conditions:
\begin{itemize}
\item The stack $X_{0}$ is a $\mathbb{U}$-small disjoint union of $(n-1)$-geometric
stacks;

\item Each of the two natural morphisms $d_{0},d_{1} : X_{1} \longrightarrow X_{0}$ is a
$(n-1)$-covering family (with one element).

\end{itemize}
\end{thm}

\textit{Sketch of proof:} Suppose that $X_{*}$ is a Segal groupoid satisfying the two conditions of the
theorem. Then, using \cite[Prop. $1.6$]{se} and corollary \ref{cex}, one checks that the natural
morphism $X_{0} \longrightarrow |X_{*}|$
is such that $X_{0}\times^{h}_{|X_{*}|}X_{0}\simeq X_{1}$. From this one deduces easily that if
$\{H_{i} \longrightarrow X_{0}\}_{i \in I}$ is an $(n-1)$-atlas for $X_{0}$, then
$\{H_{i} \longrightarrow X_{0} \longrightarrow |X_{*}|\}_{i \in I}$
is again an $(n-1)$-atlas for $|X_{*}|$. One can also check that the following diagram is homotopy cartesian
$$\xymatrix{
|X_{*}| \ar[r] & |X_{*}|\times |X_{*}| \\
X_{1} \ar[u] \ar[r] & X_{0}\times X_{0}. \ar[u]}$$
This implies that the diagonal of $|X_{*}|$ is $(n-1)$-representable and finally that $|X_{*}|$ is
$n$-geometric.

For the other implication, let $F$ be an $n$-geometric stack and $\{H_{i} \longrightarrow F\}_{i \in I}$
an $(n-1)$-atlas. Let $X_{0}:=\coprod_{i \in I}H_{i} \longrightarrow F$ be the induced morphism.
The homotopy nerve of $X_{0} \longrightarrow F$ is a simplicial object
$X_{*} \in Ho(s\mathcal{C}_{\mathbb{U}}-Aff^{\sim,\tau})$
such that $|X_{*}|\simeq F$. Furthermore, the fact that $F$ is $n$-geometric implies that
$X_{*}$ satisfies the two conditions of the theorem. \hfill $\Box$ \\

\end{subsection}

\end{section}

\begin{section}{Applications and perspectives}

In this last section, we present two applications of our theory.
The first one is an approach to $DG$-\textit{schemes} in which we
interpret them as \textit{geometric stacks over the model category
of complexes}. The second application is a definition of \'{e}tale
$K$-theory of $E_{\infty}$-ring spectra. Several other
applications will appear in a forthcoming version.

\begin{subsection}{An approach to DG-schemes}

For this paragraph let $k$ be a commutative ring with unit, $C(k)$ the symmetric
monoidal model category of complexes of $k$-modules in $\mathbb{V}$ and $C(k)_{\mathbb{U}}$ the
full sub-category of $\mathbb{C}$ of objects belonging to $\mathbb{U}$. We adopt the convention
that complexes are $\mathbb{Z}$-graded co-chain complexes (i.e differentials increase degrees).
As explained in the previous section, we will work with a fixed unital $E_{\infty}$-operad $\mathcal{O}$ in $C(k)_{\mathbb{U}}$.
For the sake of simplicity, we will assume that for each $n$, one has $\mathcal{O}(n)^{i}=0$ for $i>0$ (i.e. the operad
$\mathcal{O}$ is concentrated in non-positive degrees).

Applying definition \ref{d10},
we can consider the semi-model categories of affine stacks in $C(k)$ and in $C(k)_{\mathbb{U}}$
$$C(k)-Aff \qquad C(k)_{\mathbb{U}}-Aff.$$
In this special case, it is known that $C(k)-Aff$ and $C(k)_{\mathbb{U}}-Aff$ are actually \textit{model} categories
(see \cite{hin}). \\

Let us fix one of the standard Grothendieck topologies $\tau_{0}$
on the category of $k$-schemes (e.g. Zariski, Nisnevich,
\'{e}tale, faithfully flat, \dots). Starting from $\tau_{0}$, we
construct a topology $\tau$ on the model category
$C(k)_{\mathbb{U}}-Aff$ (Def. \ref{d4}) in the following way.
Recall that for $Spec\, A \in C(k)-Aff$, one can consider its
cohomology algebra $H^{*}(A)=\oplus H^{i}(A)$ which is in a
natural way a graded commutative $k$-algebra. The construction $A
\mapsto H^{*}(A)$ is of course functorial and therefore defines a
functor from $Ho(C(k)-Aff)^{op}$ to graded commutative
$k$-algebras.

The following definition was inspired by the work of K.
Behrend \cite{be}, where an \'{e}tale topology on differential
graded algebras is used.

\begin{df}\label{d30}
A $\mathbb{U}$-small family of morphisms in $C(k)_{\mathbb{U}}-Aff$
$$\{f_{i} : Spec\, A_{i} \longrightarrow Spec\, B\}_{i \in I}$$
is a $\tau$-covering if it satisfies the following two conditions:
\begin{itemize}
\item The induced family of morphisms of (usual) affine schemes
$$\{f_{i} : Spec\, H^{0}(A_{i}) \longrightarrow Spec\, H^{0}(B)\}_{i \in I}$$
is a $\tau_{0}$-covering;

\item For any $i \in I$, the induced morphism
$$H^{*}(B)\otimes_{H^{0}(B)}H^{0}(A_{i}) \longrightarrow H^{*}(A_{i})$$
is an isomorphism.

\end{itemize}

The topology $\tau$ on $C(k)_{\mathbb{U}}-Aff$, associated to the Grothendieck  topology $\tau_{0}$ on $k$-schemes, will
be called the \emph{strong} $\tau_{0}$-\emph{topology}. Covering families in $(C(k)_{\mathbb{U}}-Aff,\tau)$ will be
called \emph{strongly} $\tau_{0}$-\emph{covering families}.

\end{df}

The above definition allows one to introduce the strong Zariski (resp.
Nisnevich, \'{e}tale, faithfully flat and quasi-compact, \dots)
topology on $C(k)_{\mathbb{U}}-Aff$. The corresponding model site
will be denoted by $(C(k)_{\mathbb{U}}-Aff,\textrm{Zar})$ (resp.
$(C(k)_{\mathbb{U}}-Aff,\textrm{Nis})$, resp. $(C(k)_{\mathbb{U}}-Aff,\textrm{\'{e}t})$,
resp. $(C(k)_{\mathbb{U}}-Aff,\textrm{ffqc})$, \dots). The
associated model categories of stacks will be naturally denoted
by
$$C(k)_{\mathbb{U}}-Aff^{\sim,\textrm{Zar}} \qquad C(k)_{\mathbb{U}}-Aff^{\sim,\textrm{Nis}} \qquad C(k)_{\mathbb{U}}-Aff^{\sim,\textrm{\'{e}t}}
\qquad C(k)_{\mathbb{U}}-Aff^{\sim,\textrm{ffqc}}$$

\noindent and so on

\begin{prop}\label{p30}
For any Grothendieck topology $\tau_{0}$ on $k-Sch$ which is coarser than the faithfully flat and quasi-compact topology,
the induced strong $\tau_{0}$-topology $\tau$ on $C(k)_{\mathbb{U}}-Aff$ is sub-canonical.
\end{prop}

\textit{Proof:} Using lemma \ref{l7}, it is enough to show that
for any $\tau$-hypercover $Spec\, B_{*} \longrightarrow Spec\, A$ in $C(k)_{\mathbb{U}}-Aff$, the
natural morphism
$$A \longrightarrow holim_{[n] \in \Delta}B_{n}$$
is an equivalence of $E_{\infty}$-algebras. As the forgetful functor from the category of
$E_{\infty}$-algebras to the category of complexes commutes with homotopy limits, it is enough
to show that $A \longrightarrow holim_{[n] \in \Delta}B_{n}$ is a quasi-isomorphism of complexes of
$k$-modules. Furthermore, in the model category $C(k)$ the homotopy limits along $\Delta$ can be computed
using total complexes and therefore it is enough to show that the natural morphism of complexes
$$A \longrightarrow Tot(B_{*})$$
is a quasi-isomorphism. To prove this, we use the spectral sequence computing the cohomology of
a total complex as described e.g. in \cite[5.6]{we},
$$E_{2}^{p,q}=H^{p}(H^{q}(B_{*})) \Rightarrow H^{p+q}(Tot(B_{*})),$$
where $H^{q}(B_{*})$ is the normalized complex associated to the co-simplicial $k$-module
$([n] \mapsto H^{q}(B_{n}))$. Now, by definition of a $\tau$-hypercover and by the hypothesis on $\tau_{0}$, one
can use the $Tor$ spectral sequence (see \cite[thm. $V.7.3$]{km}) to prove
that the co-simplicial algebra $([n] \mapsto H^{*}(B_{n}))$ corresponds to a faithfully flat hypercover
of affine schemes
$$Spec\, H^{*}(B_{*}) \longrightarrow Spec\, H^{*}(A).$$
By the usual faithfully flat descent (see \cite[\S $I$]{mi}), the above spectral sequence degenerates and satisfies
$$E_{2}^{p,q}=0 \; \textrm{for} \; p\neq 0, \qquad E_{2}^{0,q}=H^{q}(A).$$
This in turns implies that $A \longrightarrow Tot(B_{*})$ is a quasi-isomorphism. \hfill $\Box$ \\

\bigskip
\bigskip

We recall from \cite{ck1} the notion of $DG$-scheme. We will actually adopt a slightly different definition which
is adapted to the case of an arbitrary base ring $k$. In the case $k$ is a field of characteristic zero, our notion and that of \cite{ck1} are homotopically equivalent (see below).  \\

Let $X$ be a $k$-scheme (all schemes will be separated and quasi-compact) and $CQCoh(\mathcal{O}_{X})$
its category of complexes of quasi-coherent $\mathcal{O}_{X}$-modules. This category is
an algebra over the symmetric monoidal category $C(k)$, therefore it makes sense
to talk about $E_{\infty}$-algebras in $CQCoh(\mathcal{O}_{X})$ (see \cite{sp}).

\begin{df}\label{d31}
A (separated and quasi-compact) $DG$-\emph{scheme} is a pair $(X,A_{X})$ where
$X$ is a (separated and quasi-compact) $k$-scheme and $A_{X}$ is a $E_{\infty}$-algebra
in $CQCoh(\mathcal{O}_{X})$ satisfying the following two conditions:
\begin{itemize}
\item
$A_{X}$ is concentrated in non-positive degrees (i.e.  $A_{X}^{i}=o$ for $i>0$);

\item The unit morphism $\mathcal{O}_{X} \longrightarrow A_{X}^{0}$ is
an isomorphism.

\end{itemize}

A \emph{morphism} between $DG$-schemes  $f : (X,A_{X}) \longrightarrow (Y,A_{Y})$ is
the data of a morphism of schemes $f : X \longrightarrow Y$ together with
a morphism of $E_{\infty}$-algebras in $CQCoh(\mathcal{O}_{X})$,
$f^{*}(A_{Y}) \longrightarrow A_{X}$.
\end{df}

For a $DG$-scheme $(X,A_{X})$, the cohomology sheaf $H^{0}(A_{X})$ is a quasi-coherent
$\mathcal{O}_{X}$-algebra whose associated $X$-affine scheme will be denoted
by
$$\textrm{H}^{0}(X,A_{X}):=Spec\, H^{0}(A_{X}) \longrightarrow X.$$
Actually, as $A_{X}^{0}\simeq \mathcal{O}_{X}$ and $A_{X}^{1}=0$, the scheme $\textrm{H}^{0}(X,A_{X})$ is a closed
sub-scheme of $X$.
The cohomology sheaves $H^{*}(A_{X})$ are naturally quasi-coherent $H^{0}(A_{X})$-modules and therefore
correspond to quasi-coherent sheaves on the sub-scheme $ \textrm{H}^{0}(X,A_{X})$. They will still be denoted
by $H^{*}(A_{X})$.

\begin{df}\label{d32}
A morphism of $DG$-schemes $f : (X,A_{X}) \longrightarrow (Y,A_{Y})$
is a \emph{quasi-isomorphism} if it satisfies the following two conditions:
\begin{itemize}
\item The induced morphism of schemes $\mathrm{H}^{0}(f) : \mathrm{H}^{0}(X,A_{X}) \longrightarrow \mathrm{H}^{0}(Y,A_{Y})$
is an isomorphism;

\item The natural morphism of quasi-coherent sheaves on $\mathrm{H}^{0}(X,A_{X})\simeq \mathrm{H}^{0}(Y,A_{Y})$
$$H^{*}(A_{Y}) \longrightarrow H^{*}(A_{X})$$
is an isomorphism.

\end{itemize}

The homotopy category of $DG$-schemes is the category obtained from the category of
$DG$-schemes belongings to $\mathbb{U}$ by formally inverting the quasi-isomorphisms.
It will be denoted by $Ho(DG-Sch)$.
\end{df}

\bigskip
\bigskip

\textit{Remarks:} \begin{itemize}
\item
The category of $DG$-schemes in $\mathbb{U}$ is a $\mathbb{V}$-small category. Therefore, $Ho(DG-Sch)$
is also a $\mathbb{V}$-small category but it is not clear a priory that it is
a $\mathbb{U}$-small category.

\item
When $k$ is a field of characteristic zero, the definition of $DG$-scheme given in \cite{ck1}
is not strictly equivalent to \ref{d31}. However, it is well known that in this case the homotopy theory
of commutative differential graded algebras is equivalent to the homotopy theory
of $E_{\infty}$-algebras. This fact implies easily that
the homotopy category of $DG$-schemes as defined in \cite{ck1} (and called by the authors, the right derived category of $k$-schemes) is equivalent to
our $Ho(DG-Sch)$.

\item Let $A$ be a $E_{\infty}$-algebra in $\mathbb{U}$ such that $A^{i}=0$ for $i>0$.
As the operad $\mathcal{O}$ is concentrated in non-positive degrees, the $k$-module
$A^{0}$ carries an induced $E_{\infty}$-algebra structure. As it is a complex concentrated in degree zero, this
is then nothing else than a commutative and unital algebra structure. Moreover,  it is clear that the natural
morphism of complexes $A^{0} \longrightarrow A$ is a morphism of $E_{\infty}$-algebras. In particular, $A$
is naturally a complex of $A^{0}$-modules.
This implies
that for any $E_{\infty}$-algebra $A$ such that $A^{i}=0$ for $i>0$, one can define a $DG$-scheme
$X:=\underline{Spec}\, A$, whose underlying scheme is $Spec\, A^{0}$ and with
$A_{X}:=\widetilde{A} \in QCoh(X)$. It is clear that any $DG$-scheme $(X,A_{X})$ such that $X$ is an affine
scheme is of the form $\underline{Spec}\, A$ for some $E_{\infty}$-algebra in non-positive degrees $A$
(in fact, one has $A\simeq \Gamma(X,A_{X})$).

\end{itemize}

\begin{prop}\label{p31}
There exists a functor
$$\Theta : Ho(DG-Sch) \longrightarrow Ho(C(k)_{\mathbb{U}}-Aff^{\sim,\textrm{ffqc}})$$
such that, for any $E_{\infty}$-algebra in non-positive degrees $A$, one has
$$\Theta(\underline{Spec}\, A)\simeq \mathbb{R}Spec\, A.$$
Moreover, for every $DG$-scheme $(X,A_{X})$, the stack $\Theta(X,A_{X})$ is $1$-geometric.
\end{prop}

\textit{Sketch of Proof:} Let $(X,A_{X})$ be a $DG$-scheme and let $\{U_{i}\}_{i \in I}$ be a finite Zariski
covering of $X$ by affine schemes. Taking the nerve of this covering yields a simplicial diagram
of affine schemes
$$[n] \mapsto  \coprod_{i_{0},\dots,i_{n}\in I^{n+1}}U_{i_{0},\dots,i_{n}},$$
where $U_{i_{0},\dots,i_{n}}:=U_{i_{0}}\cap \dots \cap U_{i_{n}}$.
By restricting $A_{X}$ on each $U_{i_{0},\dots,i_{n}}$, one
actually obtains a simplicial diagram of $DG$-schemes
$$[n] \mapsto  (\coprod_{i_{0},\dots,i_{n}\in I^{n+1}}U_{i_{0},\dots,i_{n}},\coprod A_{U_{i_{0},\dots,i_{n}}}).$$
Moreover, as each $\coprod_{i_{0},\dots,i_{n}\in
I^{n+1}}U_{i_{0},\dots,i_{n}}$ is an affine scheme, this diagram
is actually the image by $\underline{Spec}$ of a co-simplicial
diagram of $E_{\infty}$-algebras or, equivalently, of a simplicial
diagram in $C(k)_{\mathbb{U}}-Aff$

$$\xymatrix @R=.3pt{F(U,X) : & \Delta^{op} \ar [r] & C(k)_{\mathbb{U}}-Aff  \\
                     &  [n] \ar@{|->} [r]        & Spec\, A_{n}
                     }
$$

Considering its image by $\mathbb{R}Spec$, this diagram induces a well
defined object in
$Ho(s(C(k)_{\mathbb{U}}-Aff^{\sim,\textrm{ffqc}}))$, the homotopy
category of simplicial objects in
$C(k)_{\mathbb{U}}-Aff^{\sim,\textrm{ffqc}}$

$$\xymatrix @R=.3pt{F(U,X) : & \Delta^{op} \ar [r] & C(k)_{\mathbb{U}}-Aff^{\sim,\textrm{ffqc}}  \\
                     &  [n] \ar@{|->} [r]        & \mathbb{R}Spec\, A_{n}.
                     }
$$

We then define
$$\Theta(U,X):=hocolim_{[n] \in \Delta^{op}} \mathbb{R}Spec\, A_{n} $$
\noindent as an object in $Ho(C(k)_{\mathbb{U}}-Aff^{\sim,\textrm{ffqc}})$. With some work, it is not difficult to verify that the stack $\Theta(U,X)$ does not depend on the choice of the affine
covering $\{U_{i}\}_{i \in I}$ and that $(X,A_{X}) \mapsto \Theta(U,X)$ defines a functor
$$\Theta : Ho(DG-Sch) \longrightarrow Ho(C(k)_{\mathbb{U}}-Aff^{\sim,\textrm{ffqc}}).$$
By construction, it is clear that $\Theta(\underline{Spec}\, A)\simeq \mathbb{R}Spec\, A$.

Finally, to prove that $\Theta(X,A_{X})$ is a $1$-geometric stack, one applies the criterion
\ref{t5}. The conditions of \ref{t5} are satisfied because by construction $\Theta(X,A_{X})$ is the geometric realization
of the Segal groupoid $[n] \mapsto \mathbb{R}Spec\, A_{n}$, for which the natural morphisms $\mathbb{R}Spec\, A_{1}
\longrightarrow \mathbb{R}Spec\, A_{0}$
are clearly strong Zariski coverings and a fortiori coverings in $C(k)_{\mathbb{U}}-Aff^{\sim,\textrm{ffqc}}$. \hfill $\Box$ \\

\bigskip

We make the following

\medskip

\begin{conj}\label{c}
The functor $\Theta$ of proposition \ref{p31} is fully faithful.
\end{conj}

This conjecture says that the homotopy theory of $DG$-schemes
can be embedded into the homotopy theory of geometric stacks over
the model category of complexes. In other words, the theory of
$DG$-schemes should be \textit{a part of algebraic geometry over the model
category of complexes}. We propose the model category of stacks
$C(k)_{\mathbb{U}}-Aff^{\sim,\textrm{ffqc}}$ as a natural setting for the
theory of $DG$-schemes and more generally, for the theory of
$DG$-stacks. One of the reasons why we believe this is a natural candidate is that in this way $DG$-schemes would appear naturally as a part of a fully-fledged homotopy theory, in the abstract modern sense of Quillen model categories. Instead, trying to obtain in a complete elementary way a homotopy structure out of
usual DG-schemes (e.g., defining the weaker structure of a category with fibrations and equivalences, by declaring smooth maps to be fibrations and quasi-isomorphisms to be equivalences, as it seems to be suggested in \cite{ka}) seems to run into difficulties and moreover it is not a priori clear what kind of flexibility such a construction could have.
\end{subsection}

\begin{subsection}{\'{E}tale $K$-theory}

The problem of definition \'{e}tale $K$-theory was raised by P.A. Ostv\ae r and we give below a possible answer. We were very delighted by the question since it looked as a particularly good test of applicability of our theory. \\

Let $Sp^{\Sigma}$ be the model category of symmetric spectra in
$\mathbb{V}$ and $Sp^{\Sigma}_{\mathbb{U}}$ its sub-model category
of objects in $\mathbb{U}$ (see \cite{hss}). The wedge product of
symmetric spectra makes $Sp^{\Sigma}$ and
$Sp^{\Sigma}_{\mathbb{U}}$ into symmetric monoidal model
categories. Applying definition \ref{d10}, we may consider the
semi-model categories $Sp^{\Sigma}_{\mathbb{U}}-Aff$ of affine
stacks over $Sp^{\Sigma}_{\mathbb{U}}$.
Again, it is know that $Sp^{\Sigma}_{\mathbb{U}}-Aff$ is actually a model category.

For each object $Spec\, A \in Sp^{\Sigma}_{\mathbb{U}}-Aff$, one can consider
the category of $A$-modules in $Sp^{\Sigma}_{\mathbb{U}}$, $Mod(A)_{\mathbb{U}}$ as defined in the previous section.
As $Sp^{\Sigma}$ satisfies the monoid axiom,
$Mod(A)_{\mathbb{U}}$ is actually a model category (with fibrations and equivalences defined
on the underlying objects) which is moreover Quillen equivalent to $Mod(QA')_{\mathbb{U}}$,
where $QA'$ is a cofibrant replacement of $A$. Therefore, in theorem \ref{t4} one does not need
to ask $A$ to be a cofibrant object in order to get a good theory of modules.

Recall from \cite[Prop. $9.10$]{sp} that the homotopy category
$Ho(Mod(A)_{\mathbb{U}})$ is a closed symmetric monoidal category.
One can therefore define the notion of strongly dualizable objects
in $Ho(Mod(A)_{\mathbb{U}})$ (following \cite[\S $III.7$]{ekmm}).
The full sub-category of $Mod(A)_{\mathbb{U}}^{c}$ consisting of
strongly dualizable objects will be denoted by
$Mod(A)_{\mathbb{U}}^{sd}$, and will be equipped with the induced
notion of cofibrations and equivalences coming from
$Mod(A)_{\mathbb{U}}$. It is not difficult to check that with this
structure, $Mod(A)_{\mathbb{U}}^{sd}$ is then a Waldhausen
category (see \cite[\S $VI$]{ekmm}). Furthermore, if $A
\longrightarrow B$ is a morphism of $E_{\infty}$-algebras in
$Sp^{\Sigma}$, then the base change functor
$$f^{*} : Mod(A)_{\mathbb{U}}^{sd} \longrightarrow Mod(B)_{\mathbb{U}}^{sd},$$
being the restriction of a left Quillen functor, preserves equivalences and cofibrations.
This makes the lax functor
$$\begin{array}{cccc}
Mod(-)_{\mathbb{U}}^{sd} : & Sp^{\Sigma}_{\mathbb{U}} & \longrightarrow & Cat \\
& Spec\, A & \mapsto & Mod(A)_{\mathbb{U}}^{sd} \\
& (f : A \rightarrow B) & \mapsto & f^{*}
\end{array}$$
into a lax presheaf of Waldhausen $\mathbb{V}$-small categories.
Applying standard strictification techniques
we deduce a presheaf of $\mathbb{V}$-simplicial sets
of $K$-theory
$$\begin{array}{cccc}
K(-) : & Sp^{\Sigma}_{\mathbb{U}} & \longrightarrow & SSet \\
& Spec\, A & \mapsto & K(Mod(A)_{\mathbb{U}}^{sd}).
\end{array}$$

\begin{df}\label{d34}
The previous presheaf will be considered as an object in $Sp^{\Sigma}_{\mathbb{U}}-Aff^{\wedge}$ and will be
called the \emph{presheaf of} $K$-\emph{theory over the symmetric monoidal model category} $Sp^{\Sigma}_{\mathbb{U}}$.
For any $Spec\, A \in Sp^{\Sigma}_{\mathbb{U}}-Aff$, we will write
$$\mathbb{K}(A):=K(Spec\, A).$$
\end{df}

\bigskip
\bigskip

\textit{Remark.} The same construction as above works if one replaces $Sp^{\Sigma}_{\mathbb{U}}$ by
a general symmetric monoidal model category allowing therefore to define the spectrum
$\mathbb{K}(A)$ for any $E_{\infty}$-algebra $A$ in a general symmetric monoidal model category. It could be interesting to look at this construction for the \textit{motivic} categories considered in \cite[14.8]{sp}.\\

\begin{df}\label{d35}
Let $\tau$ be a topology on the model category $Sp^{\Sigma}_{\mathbb{U}}-Aff$ and
$Sp^{\Sigma}_{\mathbb{U}}-Aff^{\sim,\tau}$ the associated model category of stacks.
Let $K \longrightarrow K_{\tau}$ be a fibrant replacement of $K$ in
$Sp^{\Sigma}_{\mathbb{U}}-Aff^{\sim,\tau}$.

The $K_{\tau}$-\emph{theory space} of an $E_{\infty}$-algebra $A$ in $Sp^{\Sigma}_{\mathbb{U}}$ is defined
by
$$\mathbb{K}_{\tau}(A):=K_{\tau}(Spec\, A).$$
The natural morphism $K \longrightarrow K_{\tau}$ induces a natural augmentation (localization morphism)
$$\mathbb{K}(A) \longrightarrow \mathbb{K}_{\tau}(A).$$
\end{df}

\bigskip

\textit{Remark.} Note that we have
$$\mathbb{K}_{\tau}(A)\simeq \mathbb{R}\underline{Hom}_{w,\tau}(h_{Spec\, A},K)
\simeq \mathbb{R}\underline{Hom}_{w,\tau}(\mathbb{R}Spec\, A,K).$$

\bigskip

\noindent \textbf{An application: \'{e}tale $K$-theory of $E_{\infty}$-ring spectra.} \\
One defines an \'{e}tale topology on
$Sp^{\Sigma}_{\mathbb{U}}-Aff$ by stating that a family $\{f_{i}
: Spec\, B_{i} \longrightarrow Spec\, A\}_{i \in I}$ is an
\'{e}tale covering if it satisfies the following three conditions:

\begin{enumerate}
\item For all $i \in I$, the morphism $A \longrightarrow B_{i}$ is a formally \'{e}tale morphism of
$E_{\infty}$-ring spectra (in the sense that the corresponding co-tangent complex $L_{B_{i}/A}$ of \cite[$7$]{hin}
vanishes);

\item For all $i \in I$, the $A$-algebra $B_{i}$ is finitely presented (in any reasonable sense, see e.g. \cite{ma-re} or \cite[p. 7]{ro} in the "absolute" case, for connective, $p$-complete spectra)\footnote{A precise definition would  need more polishing and insight in the general case; we expect to give all details in a forthcoming version of this paper. However we are convinced that any topologically natural definition should work well to finally give a topology.};

\item The family of base change functors
$$\{\mathbb{L}f_{i}^{*} : Ho(Mod(A)_{\mathbb{U}}) \longrightarrow Ho(Mod(B_{i})_{\mathbb{U}})\}_{i \in I}$$
is conservative i.e. a morphism in $Ho(Mod(A)_{\mathbb{U}})$ is an isomorphism if and only if,
for any $i \in I$, its image in $Ho(Mod(B_{i})_{\mathbb{U}})$ is an isomorphism.

\end{enumerate}

One can check that these conditions actually define a topology
\textit{\'{e}t} on $Sp^{\Sigma}_{\mathbb{U}}-Aff$ Therefore, using
definition \ref{d35}, one can associate to any $E_{\infty}$-ring
spectrum $A$ in $Sp^{\Sigma}_{\mathbb{U}}$ its \textit{\'{e}tale}
$K$-\textit{theory space} $\mathbb{K}_{et}(A)$.

\end{subsection}

\end{section}

\end{document}